\documentclass[a4paper,12pt]{elsarticle}
\journal{ArXiv}
\usepackage{amsmath}
\usepackage[english]{babel}
\usepackage{graphicx}
\usepackage[usenames,dvipsnames]{color} 
\usepackage[numbers]{natbib}
\usepackage[footnotesize,bf,sf,center]{caption}
\usepackage{float}
\usepackage[dvips]{epsfig}
\usepackage[bottom=2.5cm,top=2.5cm,left=3cm,right=2cm]{geometry}
\usepackage{calc}
\usepackage{soul}
\usepackage{amsmath,amsfonts,amssymb}
\usepackage{subfigure}
\usepackage{setspace}
\usepackage{amssymb}
\usepackage{mathrsfs}
\usepackage{algorithm}
\usepackage{algpseudocode}
\usepackage{dsfont}
\usepackage{epstopdf}
\usepackage{subfig}%
\usepackage[utf8]{inputenc}
\usepackage{amsfonts,amsmath}
\usepackage{nicefrac,xfrac}
\newcommand{\ve}[1]{\mbox{\boldmath $#1$}}
\newtheorem{theorem}{Theorem}[section]

\newtheorem{lemma}[theorem]{Lemma}

\newtheorem{remark}{Remark}[section]

\newdefinition{rmk}{Remark}
\newcommand{\proof} [1]{ \noindent {\bf Proof.} #1 \hfill\rule{0.5em}{1.2ex} \par\medskip}
\usepackage{graphicx}
\begin{document}
\newcommand{\red}[1]{\textcolor{red}{#1}}
\newcommand{\blue}[1]{\textcolor{blue}{#1}}
\newcommand{\Orient}[1]{\colorbox{yellow}{\bf \ ... }\textcolor{ForestGreen}{#1}\colorbox{yellow}{\bf ... \ }}
\newcommand{\vol}{\mathop{\ooalign{\hfil$V$\hfil\cr\kern0.08em--\hfil\cr}}\nolimits}
\begin{frontmatter}

  \renewcommand\arraystretch{1.0}

    \title{\textbf{Implicit-explicit schemes for incompressible flow problems with variable viscosity}}
      
    \author{
   {\bf Gabriel R.~Barrenechea} $^1$, \ 
   {\bf Ernesto Castillo} $^2$\\
        and \  
   {\bf Douglas R.~Q.~Pacheco} $^{3,}$
   \footnote[0]{{\sf Email address:} {\tt douglas.r.q.pacheco@ntnu.no}, 
     {\sf corresponding author}}\\
   {\small ${}^{1}$ Department of Mathematics and Statistics, University of Strathclyde, Glasgow, Scotland} \\
   {\small ${}^{2}$ Department of Mechanical Engineering, University of Santiago de Chile, Santiago, Chile} \\
   {\small ${}^{3}$ Department of Mathematical Sciences, NTNU, Trondheim, Norway}}


\begin{abstract}
In this work we study different Implicit-Explicit (IMEX) schemes for incompressible flow problems with variable viscosity. Unlike most previous work on IMEX schemes, which focuses on the convective part, we here focus on treating parts of the diffusive term explicitly to reduce the coupling between the velocity components. We present different, both monolithic and fractional-step, IMEX alternatives for the variable-viscosity Navier--Stokes system, analysing their theoretical and algorithmic properties.  Stability results are proven for all the methods presented, with all these results being unconditional, except for one of the discretisations using a fractional-step scheme, where a CFL condition (in terms of the problem data) is required for showing stability. Our analysis is supported by a series of numerical experiments.
\end{abstract}
\end{frontmatter}

\section{Introduction}

Variable-viscosity flow problems are relevant in many physical and technological processes. Viscosity variations can be produced, for example, by temperature or pressure gradients, by non-Newtonian rheological behaviour, or by the interaction of multiple fluid phases. 
Numerically, non-constant viscosity often leads to extra non-linearities and ill-conditioning in the equations, which can affect the convergence and the performance of linear and nonlinear solvers \cite{Carey1989,Schussnig2021JCP,Pacheco2021CMAME}. These issues get amplified in the time-dependent situation, where those solvers are used multiple (even thousands) of times as time advances. As a consequence, the last few years have seen an expansion in the literature on numerical methods for variable-viscosity incompressible flow problems \cite{Schussnig2021JCP,Pacheco2021CMAME,Deteix2018,Plasman2020,Anaya2021,Anaya2023}.

Over the last few decades, several alternatives have been proposed to reduce the computational complexity of fluid simulations. Those date back to the seminal work \cite{Chorin67}, as well as operator-splitting schemes \cite{Glowinski} and more general projection methods, such as fractional-step schemes \cite{Deteix2018,chorin1968,temam1969,Guermond2006,BadiaRamon2007,Diaz2023}. A more recent trend, which can also be combined with fractional stepping, are implicit-explicit (IMEX) methods. These have become increasingly popular in applications ranging from waves to fluids \cite{Boscarino2013,Hochbruck2021,Burman2022,Guesmi2023,Burman2023}. 

In the context of incompressible flows, IMEX usually refers to temporal discretisations that treat convection explicitly or semi-implicitly, while keeping viscous terms fully implicit \cite{John2016}, which is the case for most of the IMEX literature for flow problems. To the best of our knowledge, only few articles so far have considered IMEX methods in the presence of variable viscosity. It was recently proposed in \cite{Stiller2020,Guesmi2023} to augment the viscous term with a grad-div stabilisation \cite{Olshanskii2009}, treated explicitly. Another alternative consists on adding and subtracting a large enough constant-coefficient viscous term on one side of the equation, and then applying the IMEX time-marching method to the equivalent equation. This is the approach used in \cite{WWZS16,WZWS20,TCS22}, showing good numerical performance, although the analysis carried out so far is limited to scalar equations with constant coefficients.

For flow problems with variable viscosity, treating the diffusive term semi-implicitly can present significant computational advantages. As a matter of fact, for those flow problems the viscous term cannot be written as the standard Laplacian. This induces a coupling of the different components of the velocity (due to the presence of the symmetric gradient), which makes the linear solvers more involved and memory- and time-consuming. Thus, the possibility of devising IMEX discretisations that make the transposed velocity gradient explicit, leading to linear problems with the same sparsity pattern as a scalar transport equation, has practical (and theoretical) appeal.

In this work we further develop an idea explored in \cite{Pacheco2021CMAME}, where an alternative formulation using the viscosity gradient was used to propose an IMEX scheme decoupling the velocity components; the results indicate an improved temporal stability of the modified formulation in comparison to ``naive'' IMEX discretisations. More precisely, we study a variety of IMEX discretisations, both monolitic and fractional-step, arising from making one part of the viscous term explicit, while keeping the ``Laplacian'' part implicit. We prove stability for all the methods studied, but the norm with respect to which the stability is proven varies according to the formulation. In fact, we show that the most ``natural'' way of splitting the viscous term leads to a weak stability, while a more careful rewriting of the term to be made explicit enhances it, which is confirmed by our numerical studies. To keep the technical details to a minimum, and focus on the effects of the time discretisation, we only analyse the semi-discretised (in time) problem, and consider only temporally first-order schemes, both for the monolithic and the fractional-step alternatives.

The rest of the manuscript is organised as follows. In Section~\ref{sec_pre} we introduce the model problem, notation and useful analytical tools. In Sections~\ref{sec_mono} and \ref{sec_incremental} we derive stability estimates for monolithic and fractional-step methods, respectively. Section \ref{sec_summary} summarises the theoretical results of the paper, and in Section \ref{sec_spatial} we briefly discuss spatial discretisation matters. We test and compare numerically the different alternatives presented in Section~\ref{sec_examples}, and we draw some final remarks in Section~\ref{sec_Conclusion}.

\section{Preliminaries}\label{sec_pre}
\subsection{Model problem}
Let us consider a finite time interval $(0,T]$ and a fluid domain $\Omega\subset\mathbb{R}^{d}$, $d=2$ or $3$, with Lipschitz boundary $\Gamma=\partial\Omega$. As a model problem, we consider the variable-viscosity incompressible Navier--Stokes equations:
\begin{align}
\partial_t\ve{u} + (\nabla\ve{u})\ve{u} - \nabla\cdot[2\nu(\ve{x},t)\nabla^{\mathrm{s}}\ve{u}] + \nabla p &= \ve{f} && \text{in} \ \ \Omega\times(0,T]=:Q\, ,\label{momentum}\\
\nabla\cdot\ve{u} &= 0  && \text{in} \ \ \Omega\times(0,T]\, , \label{incompressibility}\\
\ve{u} &= \ve{0} && \text{on} \ \ \Gamma\times(0,T]\, ,\label{DirichletBC}\\
\ve{u} &= \ve{u}_0 && \text{at} \ \ t=0\, ,
\end{align}
where the unknowns are the velocity $\ve{u}$ and the pressure $p$, while the remaining quantities are problem data. The viscosity field $\nu(\ve{x},t)$ is assumed to satisfy
\begin{equation}
    \nu(\ve{x},t) \geq  \nu_{\mathrm{min}} > 0 \ \ \text{in} \ \bar{Q}\, ,
\end{equation}
where $\nu_{\mathrm{min}}$ is a known constant that represents, e.g., the minimum viscosity of a shear-thinning fluid. To simplify notation, we shall use $\nu$ instead of $\nu(\ve{x},t)$ throughout this work (but will keep in mind its variable character).

\subsection{Alternative formulations of the viscous term}\label{sec_Weak}
In incompressible flows, the divergence-free constraint \eqref{incompressibility} leads to multiple ways of describing the viscous term. Although popular in the constant viscosity scenario, the Laplacian form $\nu\Delta\ve{u}$ is not appropriate for problems with variable viscosity. For that reason, the full stress-divergence (SD) form
\begin{equation}
\nabla\cdot\left(2\nu\nabla^{\mathrm{s}}\ve{u}\right) = \nabla\cdot(\nu\nabla\ve{u}) + \nabla\cdot(\nu\nabla^{\top}\ve{u}) 
\label{SD}
\end{equation}
is normally used when $\nabla\nu \not= \ve{0}$. Yet, other possibilities exist, such as \cite{Guesmi2023}
\begin{align}
\nabla\cdot\left(2\nu\nabla^{\mathrm{s}}\ve{u}\right) = \nabla\cdot\left(2\nu\nabla^{\mathrm{s}}\ve{u}\right) - \nabla(\nu\nabla\cdot\ve{u})\, ,
\label{gradDiv}
\end{align}
or also \cite{Anaya2021}
\begin{align}\label{Anaya}
\nabla\cdot\left(2\nu\nabla^{\mathrm{s}}\ve{u}\right) &= 2\nabla^{\mathrm{s}}\ve{u}\nabla\nu - \nu\nabla\times(\nabla\times\ve{u}) \, .
\end{align}
In this work we consider the following rewriting of the viscous term
\begin{equation}
\nabla\cdot\left(2\nu\nabla^{\mathrm{s}}\ve{u}\right) \equiv \nabla\cdot\left(\nu\nabla\ve{u}\right) + \nabla^{\top}\ve{u}\nabla\nu+ \nu\nabla(\nabla\cdot\ve{u})  \\
    = \nabla\cdot\left(\nu\nabla\ve{u}\right) + \nabla^{\top}\ve{u}\nabla\nu\, ,
\label{genLapStrong}
\end{equation}
which we denote as generalised Laplacian (GL). Although equivalent at the continuous level, the writings \eqref{SD}-\eqref{genLapStrong} can (and do) lead to dramatically different numerical performances. As mentioned in the introduction, the symmetric gradient leads to a coupling of the $d$ velocity components (especially if $\nabla\nu\not=\ve{0}$). In this work our aim is to make the parts of the viscous term responsible for this coupling explicit, and to analyse the implications of such IMEX approaches.

\subsection{Useful notation, identities and inequalities}
We consider usual notation for Hilbert and Lebesgue spaces aligned, e.g., with \cite{EG21-I}. More precisely, we denote by $L^2_0(\Omega)$ the space of $L^2(\Omega)$ functions with zero mean on $\Omega$. The inner product in $L^2(\Omega)$ is denoted by $( \cdot ,\cdot )$, and we make no distinction between the product of scalar or vector/tensor-valued functions. Moreover, we denote by $\| \cdot \|$ and $\| \cdot \|_{\infty}$ the $L^2(\Omega)$ and $L^{\infty}(\bar{Q})$ norms, respectively. In our analyses, we will assume that $\ve{f}\in [L^2(Q)]^d$ and $\ve{u}_0\in [H^1_0(\Omega)]^d$.

We will carry out a discrete-in-time analysis. Approximate or discrete values of the variables at the different instants will be defined with sub-indices: $\ve{u}_{n}$, for instance, denotes the velocity approximation at the $n$-th time step. We will often use the identity
\begin{align}
    2(\ve{u}_{n+1}-\ve{u}_{n},\ve{u}_{n+1}) = \|\ve{u}_{n+1}\|^2 - \|\ve{u}_{n}\|^2 + \|\delta\ve{u}_{n+1}\|^2 \, ,
    \label{identity}
\end{align}
where $\delta\ve{u}_{n+1}:=\ve{u}_{n+1}-\ve{u}_{n}$. 
A classical and useful fact is the skew-symmetry of the convective term. More precisely, if $\ve{v}\in [H^1_0(\Omega)]^d$ is solenoidal, there holds
\begin{align}
((\nabla\ve{w})\ve{v},\ve{w}) = 0 \ \ \text{for all $\ve{w}\in [H^1_0(\Omega)]^d$.} 
\label{convective}
\end{align}
If $\ve{v}$ is \textit{not} divergence-free, we will consider the following skew-symmetrised form $\mathbf{c}$ of the convective derivative:
\begin{align}
   \mathbf{c}(\ve{v},\ve{w}) :=(\nabla\ve{w})\ve{v} + \frac{1}{2}(\nabla\cdot\ve{v})\ve{w}\, ,
   \label{skewSymmetrisation}
\end{align}
which enjoys the property
\begin{align}
   (\mathbf{c}(\ve{v},\ve{w}),\ve{w}) = 0 \ \ \text{for all} \ \ve{w},\ve{v} \in [H^1_0(\Omega)]^d \, .
   \label{skewSymmetric}
\end{align}
For details on different forms of the convective term and their properties, we refer the reader to, e.g., \cite{John2016}.

To construct IMEX schemes, we can combine backward differentiation formulas in time with extrapolation rules of matching order of  consistency. Our analysis will be restricted to the first-order case, for which extrapolation means simply replacing a certain quantity at $t_{n+1}$ by its previous value at $t_n$, and
\begin{align*}
    \partial_t\ve{u}|_{t=t_{n+1}} \approx \frac{1}{\tau}(\ve{u}_{n+1}-\ve{u}_{n})\, ,
\end{align*}
where $\tau>0$ denotes the time-step size, which will be considered constant, for simplicity -- although all the results in this paper extend to variable $\tau$ in a straightforward manner.

We will analyse both monolithic and fractional-step schemes. The stability analysis will rely on different versions of the discrete Gronwall inequality. The proof of the two results below can be found in \cite[Lemma~5.1]{Heywood1990} (see also \cite[Lemma~A.56]{John2016}).
\begin{lemma}[Discrete Gronwall inequality]\label{Lem:Gronwall-unconditional}
Let $N\in\mathbb{N}$, and $\alpha,B,a_{n},b_{n},c_{n}$ be non-negative numbers for $n=1,\ldots,N$. Let us suppose that these numbers satisfy
\begin{align}\label{N-1}
    a_{N} + \sum_{n=1}^{N}b_n \leq B + \sum_{n=1}^{N}c_n + \alpha\sum_{n=1}^{N-1}a_n  \, .
\end{align}
Then, the following inequality holds:
\begin{align}
    a_{N} + \sum_{n=1}^{N}b_n \leq \mathrm{e}^{\alpha N}\left( B + \sum_{n=1}^{N}c_n\right) \ \ \text{for} \ \ N\geq 1 \, .
\end{align}
\end{lemma}

\begin{lemma}[Discrete Gronwall lemma, conditional version]
\label{Lem:Gronwall-conditional}
Let $N\in\mathbb{N}$, and \\
$\alpha,B,a_{n},b_{n},c_{n}$ be non-negative numbers for $n=1,\ldots,N$. Let us suppose that these numbers satisfy
\begin{align}\label{N}
    a_{N} + \sum_{n=1}^{N}b_n \leq B + \sum_{n=1}^{N}c_n + \alpha\sum_{n=1}^{N}a_n \, ,
\end{align}
with $\alpha < 1$. Then, the following inequality holds:
\begin{align}
    a_{N} + \sum_{n=1}^{N}b_n \leq \mathrm{e}^{\frac{\alpha N}{1-\alpha}}\left( B + \sum_{n=1}^{N}c_n\right)\, .
\end{align}
\end{lemma}
Notice that the only difference between \eqref{N-1} and \eqref{N} is where the last sum on the right-hand side stops ($N-1$ or $N$).

\section{Monolithic schemes}\label{sec_mono}

In this section we will study three different IMEX schemes. All of them share the property that the velocity and the pressure are computed simultaneously at each time step (hence the term ``monolithic'' to describe them). The proof of their stability will be based on the discrete Gronwall inequality given in Lemma~\ref{Lem:Gronwall-unconditional}. Hence, every time we refer to the {\sl discrete Gronwall inequality} in this section, we will be referring to Lemma~\ref{Lem:Gronwall-unconditional}.

\subsection{Implicit stress-divergence formulation}
We will start by considering the simplest possible IMEX scheme, where (other than the convective velocity) only the viscosity is made explicit. In such a case, the time-discrete momentum equation reads
\begin{equation}
\frac{1}{\tau}(\ve{u}_{n+1}-\ve{u}_{n})-\nabla\cdot\left(2\nu_{n}\nabla^{\mathrm{s}}\ve{u}_{n+1}\right) + (\nabla\ve{u}_{n+1})\ve{u}_n + \nabla p_{n+1} =  \ve{f}_{n+1}\, .
    \label{strongMomentumMono}
\end{equation}
The standard variational problem for the $(n+1)$-th time step is given by: find $(\ve{u}_{n+1},p_{n+1})\in [H^1_0(\Omega)]^d\times L^2_0(\Omega)$ such that 
\begin{equation}
\begin{split}
\frac{1}{\tau}\left(\ve{u}_{n+1}-\ve{u}_{n},\ve{w}\right)+\left(2\nu_{n}\nabla^{\mathrm{s}}\ve{u}_{n+1},\nabla^{\mathrm{s}}\ve{w}\right) + \big((\nabla\ve{u}_{n+1})\ve{u}_n,\ve{w} \big) -\left(p_{n+1},\nabla\cdot\ve{w}\right) &=  \left( \ve{f}_{n+1},\ve{w} \right), \\
\left(q,\nabla\cdot\ve{u}_{n+1}\right) &= 0\,,
\label{fullyCoupled}
\end{split}
\end{equation}
for all $(\ve{w},q)\in [H^1_0(\Omega)]^d\times L^2_0(\Omega)$. 

\begin{remark} The decision to make the viscosity explicit in \eqref{fullyCoupled} stems from the fact that in realistic situations (such as generalised Newtonian flows), the viscosity is seldom known beforehand, being instead velocity-dependent. This is why we  consider the first-order extrapolation $\nu_{n+1}\approx \nu_n$. 
\end{remark}

The scheme \eqref{fullyCoupled} is rather classical, and well-understood. Nevertheless, we now reproduce its stability result (and its proof) for completeness, and to highlight the main differences that appear in more involved schemes presented later.

\begin{lemma}[Stability of the implicit stress-divergence formulation]
Let us suppose that the time-step size is given by $\tau=T/N$, $N\ge 2$. Then, the scheme \eqref{fullyCoupled} satisfies the following stability inequality:
\begin{equation}
    \|\ve{u}_{N}\|^2 + 4\sum_{n=1}^{N}\tau\big\|\sqrt{\nu_{n-1}}\, \nabla^{\mathrm{s}}\ve{u}_{n}\big\|^2 \leq \mathrm{e}^{\varepsilon}\Bigg[\left(1+\varepsilon\frac{\tau}{T}\right)\|\ve{u}_{0}\|^2 + \left(\tau+\frac{ T}{\varepsilon}\right)\sum_{n=1}^{N}\tau\|\ve{f}_{n}\|^2\Bigg]\, ,
\label{stabilitySD}
\end{equation}
where $\varepsilon>0$ is an arbitrary constant. 
\end{lemma}
\proof{Setting $\ve{w} = 2\tau\ve{u}_{n+1}$ and $q = 2\tau p_{n+1}$ in \eqref{fullyCoupled} and adding the resulting equations yields
\begin{align*}
    2(\ve{u}_{n+1}-\ve{u}_{n},\ve{u}_{n+1}) + 4\tau\big\|\sqrt{\nu_{n}}\, \nabla^{\mathrm{s}}\ve{u}_{n+1}\big\|^2 + 2\tau((\nabla\ve{u}_{n+1})\ve{u}_n,\ve{u}_{n+1}) = 2\tau(\ve{f}_{n+1},\ve{u}_{n+1})\, .
\end{align*}
The convective term vanishes due to its skew-symmetric property \eqref{convective}. Applying \eqref{identity} and the Cauchy-Schwarz and Young inequalities, we obtain
\begin{align*}
    \|\ve{u}_{n+1}\|^2 - \|\ve{u}_{n}\|^2 + 4\tau\big\|\sqrt{\nu_{n}}\, \nabla^{\mathrm{s}}\ve{u}_{n+1}\big\|^2  &\leq 2\tau(\ve{f}_{n+1},\ve{u}_{n+1}) \\
    &=2\tau(\ve{f}_{n+1},\ve{u}_{n}+\delta\ve{u}_{n+1}) \\
    &\leq 2\tau\|\delta \ve{u}\|\,\|\ve{f}_{n+1}\|  + 2\tau\|\ve{f}_{n+1}\|\, \|\ve{u}_{n}\|  \\
    &\leq \|\delta\ve{u}_{n+1}\|^2 + \Big(\tau+\frac{T}{\varepsilon}\Big)\tau\|\ve{f}_{n+1}\|^2 + \tau\frac{\varepsilon}{T}\|\ve{u}_{n}\|^2  ,
\end{align*}
where $\varepsilon>0$ is arbitrary. Rearranging terms we get
\begin{align*}
    \|\ve{u}_{n+1}\|^2-\|\ve{u}_{n}\|^2 + 4\tau\big\|\sqrt{\nu_{n}}\, \nabla^{\mathrm{s}}\ve{u}_{n+1}\big\|^2 \leq  \varepsilon\frac{\tau}{T}\|\ve{u}_{n}\|^2 + \Big(\tau^2+\frac{\tau T}{\varepsilon }\Big)\|\ve{f}_{n+1}\|^2 \, .
\end{align*}
Adding this equality from $n=0$ to $n=N-1$ yields
\begin{align*}
    \|\ve{u}_{N}\|^2 + 4\tau\sum_{n=1}^{N}\big\| \sqrt{\nu_{n-1}}\,\nabla^{\mathrm{s}}\ve{u}_{n}\big\|^2 \leq \left(1+\varepsilon\frac{\tau}{T}\right)\|\ve{u}_{0}\|^2 + \varepsilon\frac{\tau}{T}\sum_{n=1}^{N-1}\|\ve{u}_{n}\|^2 + \Big(\tau+\frac{T}{ \varepsilon }\Big)\sum_{n=1}^{N}\tau\|\ve{f}_{n}\|^2  ,
\end{align*}
and the result follows as a direct application of the discrete Gronwall lemma.}

\begin{remark}
The above result applied to $\ve{f}=\ve{0}$ leads to the fact that \eqref{fullyCoupled} is strongly stability preserving, i.e., 
\begin{equation*}
  \|\ve{u}_{N}\|^2 + 4\sum_{n=1}^{N}\tau\big\|\sqrt{\nu_{n-1}}\, \nabla^{\mathrm{s}}\ve{u}_{n}\big\|^2 \leq \|\ve{u}_{0}\|^2   
\end{equation*}
The main question, to be addressed next, is whether and how we can preserve stability even if a part of the viscous term is treated explicitly. 
\end{remark}

\subsection{Implicit-explicit stress-divergence formulations}
The simplest, naive way to decouple the velocity components is to split the viscous term into two parts, as in \eqref{SD}, and make the second term explicit. This choice leads to the following semi-discrete scheme: 
\begin{equation}
\frac{1}{\tau}(\ve{u}_{n+1}-\ve{u}_{n})-\nabla\cdot\left(\nu_{n}\nabla\ve{u}_{n+1}\right) + (\nabla\ve{u}_{n+1})\ve{u}_n + \nabla p_{n+1} = \nabla\cdot\left(\nu_{n}\nabla^{\top}\ve{u}_{n}\right) +  \ve{f}_{n+1}\, ,
\label{SDnaive}
\end{equation}
with the corresponding weak form: find $(\ve{u}_{n+1},p_{n+1})\in [H^1_0(\Omega)]^d\times L^2_0(\Omega)$ such that
\begin{equation*}
\begin{split}
\frac{1}{\tau}(\ve{u}_{n+1}-\ve{u}_{n},\ve{w})+\left(\nu_{n}\nabla\ve{u}_{n+1},\nabla\ve{w}\right) + \big((\nabla\ve{u}_{n+1})\ve{u}_n,\ve{w} \big)-\left(p_{n+1},\nabla\cdot\ve{w}\right) &= \\
-\left(\nu_{n}\nabla^{\top}\ve{u}_{n},\nabla\ve{w}\right) + (\ve{f}_{n+1},\ve{w})\,, &\\
\left(q,\nabla\cdot\ve{u}_{n+1}\right) &= 0\,,
\end{split}
\end{equation*}
for all $(\ve{w},q)\in [H^1_0(\Omega)]^d\times L^2_0(\Omega)$.

Testing with the same functions as in the proof of the last result, we get to
\begin{align*}
     \|\ve{u}_{n+1}\|^2  + 2\tau\big\|\sqrt{\nu_n}\,\nabla\ve{u}_{n+1} \big\|^2 + \|\delta\ve{u}_{n+1}\|^2 = \|\ve{u}_{n}\|^2 - 2\tau(\nu_n\nabla^{\top}\ve{u}_{n},\nabla\ve{u}_{n+1}) + 2\tau(\ve{f}_{n+1},\ve{u}_{n+1}) \,.
\end{align*}
In this case, an extra term (that requires bounding) appears on the right-hand side:
\begin{equation*}
-2(\nu_n\nabla^{\top}\ve{u}_{n},\nabla\ve{u}_{n+1}) \leq 2\big\|\sqrt{\nu_n}\,\nabla\ve{u}_n\big\|\,\big\|\sqrt{\nu_n}\,\nabla\ve{u}_{n+1}\big\| \leq \big\|\sqrt{\nu_n}\,\nabla\ve{u}_n\big\|^2 + \big\|\sqrt{\nu_n}\,\nabla\ve{u}_{n+1}\big\|^2\,,
\end{equation*}
where we have used the Cauchy-Schwarz and Young inequalities.
Then, estimating the forcing term as done for the implicit case, we get the following bound at step $(n+1)$:
\begin{equation}\label{20-1}
     \|\ve{u}_{n+1}\|^2  + \tau\big\|\sqrt{\nu_n}\,\nabla\ve{u}_{n+1} \big\|^2  \leq \|\ve{u}_{n}\|^2 + \tau\big\|\sqrt{\nu_n}\,\nabla\ve{u}_{n}\big\|^2 + \frac{\tau}{T}\|\ve{u}_{n}\|^2 + (\tau^2 + \tau T)\|\ve{f}_{n+1}\|^2\,.
\end{equation}

At this point, there is no guarantee that the viscous term on the left-hand side will be able to control the one on the right-hand side, since the viscosity on the left-hand side term is $\nu_n$, not $\nu_{n+1}$. Some stability might still be attained if $\nu_{n+1} \leq \nu_n$, or if $\tau$ is small enough so that $\nu_{n+1}\approx\nu_{n}$; otherwise, the scheme could become unstable, especially for rapidly varying viscosity fields. This may explain the poor stability observed for a related implicit-explicit SD formulation in \cite{Pacheco2021CMAME}, see also Section \ref{sec_cavity} below.

According to the analysis above, guaranteeing stability of the IMEX-SD formulation would require knowing $\nu_{n+1}$ already at the $(n+1)$-th time-step. Although this is seldom realistic, a scheme that requires $\nu_{n+1}$ can still make practical sense if combined, for instance, with a Picard-like method linearising the current viscosity. Therefore, let us now present a scheme that assumes $\nu_{n+1}$ as already known at the point of computing $(\ve{u}_{n+1},p_{n+1})$. Consider the first-order extrapolation
\begin{align*}
   \nu_{n+1}\nabla^{\top}\ve{u}_{n+1} = \sqrt{\nu_{n+1}}\sqrt{\nu_{n+1}}\, \nabla^{\top}\ve{u}_{n+1}  \approx \sqrt{\nu_{n+1}}\sqrt{\nu_{n}}\, \nabla^{\top}\ve{u}_{n} \, ,
\end{align*}
from which we can propose the following alternative to \eqref{SDnaive}:
\begin{equation*}
\frac{1}{\tau}(\ve{u}_{n+1}-\ve{u}_{n})-\nabla\cdot\left(\nu_{n+1}\nabla\ve{u}_{n+1}\right) + (\nabla\ve{u}_{n+1})\ve{u}_n + \nabla p_{n+1} = \nabla\cdot\left(\sqrt{\nu_{n+1}}\sqrt{\nu_{n}}\, \nabla^{\top}\ve{u}_{n}\right) +  \ve{f}_{n+1}
\end{equation*}
and the corresponding weak form: find $(\ve{u}_{n+1},p_{n+1})\in [H^1_0(\Omega)]^d\times L^2_0(\Omega)$ such that
\begin{equation}
\begin{split}
\frac{1}{\tau}(\ve{u}_{n+1}-\ve{u}_{n},\ve{w})+\left(\nu_{n+1}\nabla\ve{u}_{n+1},\nabla\ve{w}\right) + ((\nabla\ve{u}_{n+1})\ve{u}_n,\ve{w} )-\left(p_{n+1},\nabla\cdot\ve{w}\right) &= \\
-\left(\sqrt{\nu_{n}}\,\nabla^{\top}\ve{u}_{n},\sqrt{\nu_{n+1}}\,\nabla\ve{w}\right) + (\ve{f}_{n+1},\ve{w})\, , & \\
(q,\nabla\cdot\ve{u}_{n+1}) &= 0\,,
\label{SDsqrt}
\end{split}
\end{equation}
for all $(\ve{w},q)\in [H^1_0(\Omega)]^d\times L^2_0(\Omega)$. Based on that, we will prove the following stability estimate.
\begin{theorem}[Stability of an IMEX stress-divergence formulation]\label{Th:Monosqrt}
Let us suppose that the time-step size is given by $\tau=T/N$, $N\ge 2$. Then, the scheme \eqref{SDsqrt} satisfies the following stability inequality
\begin{equation}
    \|\ve{u}_{N}\|^2 + \tau\big\| \sqrt{\nu_{N}}\,\nabla\ve{u}_{N}\big\|^2 + \varepsilon\sum_{n=1}^{N-1}\frac{\tau^2}{T}\big\| \sqrt{\nu_{n}}\,\nabla\ve{u}_{n}\big\|^2 \leq \mathrm{e}^{\varepsilon}\Bigg[C_0 + \Big(\tau + \frac{T}{\varepsilon}\Big)\sum_{n=1}^{N}\tau\|\ve{f}_{n}\|^2\Bigg] ,
\label{stabilityMonoSqrt}
\end{equation}
where $\varepsilon>0$ is an arbitrary constant, and 
\begin{align*}
    C_0 = \left(1+\varepsilon\frac{\tau}{T}\right)\|\ve{u}_{0}\|^2 + \tau\|\sqrt{\nu_0}\,\nabla\ve{u}_0 \|^2\, .
\end{align*}
\end{theorem}
\proof{Using the same techniques as in the proof of \eqref{20-1} we can show that
\begin{align*}
     \|\ve{u}_{n+1}\|^2  -\|\ve{u}_{n}\|^2 + \tau\big\|\sqrt{\nu_{n+1}}\,\nabla\ve{u}_{n+1} \big\|^2  \leq  \tau\big\|\sqrt{\nu_n}\,\nabla\ve{u}_{n}\big\|^2 + \varepsilon\frac{\tau}{T}\|\ve{u}_{n}\|^2 + \Big(\tau + \frac{ T}{\varepsilon }\Big)\tau\|\ve{f}_{n+1}\|^2  .
\end{align*}
Adding up from $n=0$ to $n=N-1$ gives 
\begin{align}
    \|\ve{u}_{N}\|^2 + \tau\big\| \sqrt{\nu_{N}}\,\nabla\ve{u}_{N}\big\|^2  \leq \nonumber\\ \left(1+\varepsilon\frac{\tau}{T}\right)\|\ve{u}_{0}\|^2 + \tau\big\|\sqrt{\nu_0}\,\nabla\ve{u}_0\big\|^2 + \varepsilon\frac{\tau}{T}\sum_{n=1}^{N-1}\|\ve{u}_n\|^2 + \Big(\tau + \frac{T}{\varepsilon}\Big)\sum_{n=1}^{N}\tau\|\ve{f}_{n}\|^2\,.\label{Sec3.2.last}
\end{align}
The proof is completed by adding
\begin{align*}
    \varepsilon\frac{\tau}{T}\sum_{n=1}^{N-1}\tau\big\| \sqrt{\nu_{n}}\,\nabla\ve{u}_{n}\big\|^2
\end{align*}
to both sides of \eqref{Sec3.2.last} and using the discrete Gronwall inequality with $a_n=\|\ve{u}_{N}\|^2 + \tau\big\| \sqrt{\nu_{N}}\,\nabla\ve{u}_{N}\big\|^2$, $b_n= \varepsilon\frac{\tau}{T}\tau\big\| \sqrt{\nu_{N}}\,\nabla\ve{u}_{N}\big\|^2$, and $c_n$ gathering the terms involving $\ve{u}_0$ and $\ve{f}$ on the right-hand side of \eqref{Sec3.2.last}.}

\subsection{Implicit-explicit generalised Laplacian formulation}
Due to the squared time-step size accompanying the viscous sum in \eqref{stabilityMonoSqrt}, the stability estimate proven in Theorem~\ref{Th:Monosqrt} is somewhat weaker than the one proved for \eqref{fullyCoupled}. It is, still, an unconditional stability attained \textit{despite} the explicit treatment of the coupling term. Of course, this comes at the cost of requiring a-priori knowledge of the viscosity field. To circumvent this requirement, we now present a different IMEX scheme based on replacing the viscous term by the generalised Laplacian version \eqref{genLapStrong}. In that case, treating the coupling term explicitly will lead to the semi-discrete momentum equation
\begin{equation*}
\frac{1}{\tau}(\ve{u}_{n+1}-\ve{u}_{n})-\nabla\cdot\left(\nu_{n}\nabla\ve{u}_{n+1}\right) + (\nabla\ve{u}_{n+1})\ve{u}_n + \nabla p_{n+1} = \nabla^{\top}\ve{u}_{n}\nabla\nu_{n} +  \ve{f}_{n+1}\, .
\end{equation*}
The corresponding variational problem is: find $(\ve{u}_{n+1},p_{n+1})\in [H^1_0(\Omega)]^d\times L^2_0(\Omega)$ such that
\begin{equation}
\begin{split}
\frac{1}{\tau}(\ve{u}_{n+1}-\ve{u}_{n},\ve{w})+(\nu_{n}\nabla\ve{u}_{n+1},\nabla\ve{w}) +((\nabla\ve{u}_{n+1})\ve{u}_{n},\ve{w}) -(p_{n+1},\nabla\cdot\ve{w}) &=  \\
(\nabla^{\top}\ve{u}_{n}\nabla\nu_n,\ve{w}) + (\ve{f}_{n+1},\ve{w})\, , & \\
(q,\nabla\cdot\ve{u}_{n+1}) &= 0\,,
\label{decoupledGL}
\end{split}
\end{equation}
for all $(\ve{w},q)\in [H^1_0(\Omega)]^d\times L^2_0(\Omega)$.  The stability of \eqref{decoupledGL} is stated next.

\begin{theorem}[Stability of the IMEX generalised Laplacian formulation]
Let us suppose that $\nabla\nu\in [L^\infty(\bar{Q})]^d$, and that the time-step size is given by $\tau=T/N$ with $N\ge 2$.
Then,  the IMEX scheme \eqref{decoupledGL} satisfies the following stability inequality
    \begin{equation}
    \|\ve{u}_{N}\|^2 + 
    \sum_{n=1}^{N}\tau\big\|\sqrt{\nu_{n-1}}\,\nabla\ve{u}_{n}\big\|^2 \leq  M\Bigg[C_0  + \left(2\tau + \frac{T}{\varepsilon}\right)\sum_{n=1}^{N}\tau\|\ve{f}_n\|^2\Bigg], \label{mainResult}
    \end{equation}
    where
   \begin{align*}
       M =\mathrm{exp}\Big(\varepsilon + \frac{2T\|\nabla\nu \|_{\infty}^2}{\nu_{\mathrm{min}}}\Big)\, ,\quad 
       C_0 = \Big(1+\frac{\varepsilon\tau}{ T}+\frac{2\|\nabla\nu \|^2_{\infty}\tau}{\nu_{\mathrm{min}}}\Big)(\|\ve{u}_0 \|^2 +\tau\nu_{\mathrm{min}}\|\nabla\ve{u}_0\|^2)\, ,
   \end{align*} 
and $\varepsilon>0$ is an arbitrary constant.
\end{theorem}
\proof{Set $\ve{w}=2\tau\ve{u}_{n+1}$ and $q = 2\tau p_{n+1}$ in \eqref{decoupledGL} to get
\begin{align}
     \|\ve{u}_{n+1}\|^2   - \|\ve{u}_{n}\|^2+ \|\delta\ve{u}_{n+1}\|^2 + 2\tau\big\|\sqrt{\nu_n}\,\nabla\ve{u}_{n+1}\big\|^2 &= 2\tau(\nabla^{\top}\ve{u}_{n}\nabla\nu_n,\ve{u}_{n+1}) + 2\tau(\ve{f}_{n+1},\ve{u}_{n+1})\,.\label{CC}
\end{align} 
We now bound the terms on the right-hand side of the inequality above. Using the Cauchy-Schwarz, H\"older, and Young inequalities, we get 
\begin{align}
2\tau(\nabla^{\top}\ve{u}_{n}\nabla\nu_n,\ve{u}_{n+1}) &=2\tau(\nabla^{\top}\ve{u}_{n}\nabla\nu_n,\delta\ve{u}_{n+1}+\ve{u}_{n})\nonumber\\
 &\leq 2\tau\|\delta\ve{u}_{n+1} \|\,\|\nabla\nu\|_{\infty} \|\nabla\ve{u}_n \|  + 2\tau\|\nabla\nu\|_{\infty} \|\nabla\ve{u}_n \|\, \|\ve{u}_{n} \| \nonumber\\
 &\leq \frac{1}{2}\|\delta\ve{u}_{n+1} \|^2 + 2\tau^2\|\nabla\nu\|_{\infty}^2 \|\nabla\ve{u}_n \|^2 + \tau\nu_{\mathrm{min}} \|\nabla\ve{u}_n \|^2  + \frac{\tau\|\nabla\nu\|_{\infty}^2}{\nu_{\mathrm{min}}} \|\ve{u}_{n}\|^2\,, \label{AA}
\end{align}
and 
\begin{align}
    2\tau(\ve{f}_{n+1},\ve{u}_{n+1}) &=2\tau(\ve{f}_{n+1},\ve{u}_{n}+\delta\ve{u}_{n+1}) \nonumber\\
    &\leq 2\tau\|\delta\ve{u}_{n+1}\|\,\|\ve{f}_{n+1}\|  + 2\tau\|\ve{f}_{n+1}\|\, \|\ve{u}_{n}\|  \nonumber\\
    &\leq \frac{1}{2}\|\delta\ve{u}_{n+1}\|^2 + \left(2\tau+\frac{T}{\varepsilon }\right)\tau\|\ve{f}_{n+1}\|^2 + \frac{\varepsilon \tau}{T}\|\ve{u}_{n}\|^2\,.\label{BB}
\end{align}
Inserting \eqref{AA} and \eqref{BB} in \eqref{CC} gives
\begin{align*}
     &\|\ve{u}_{n+1}\|^2 - \|\ve{u}_{n}\|^2 + 2\tau\big\|\sqrt{\nu_{n}}\,\nabla\ve{u}_{n+1}\big\|^2 \\ 
     &\leq \tau\nu_{\mathrm{min}}\|\nabla\ve{u}_{n}\|^2+ \Big(\frac{\varepsilon\tau}{ T}+\frac{\tau\|\nabla\nu\|_{\infty}^2}{\nu_{\mathrm{min}}}\Big)\|\ve{u}_{n}\|^2 + \frac{2\tau\|\nabla\nu\|_{\infty}^2}{\nu_{\mathrm{min}}}\tau\nu_{\mathrm{min}}\|\nabla\ve{u}_{n}\|^2 + \left(2\tau+\frac{T}{\varepsilon}\right)\tau\|\ve{f}_{n+1}\|^2\\
     &\leq \tau\nu_{\mathrm{min}}\|\nabla\ve{u}_{n}\|^2+ \Big(\frac{\varepsilon\tau}{T}+\frac{2\tau\|\nabla\nu\|_{\infty}^2}{\tau\nu_{\mathrm{min}}}\Big)(\|\ve{u}_{n}\|^2+\tau\nu_{\mathrm{min}}\|\nabla\ve{u}_{n}\|^2)  + \left(2\tau+\frac{T}{\varepsilon}\right)\tau\|\ve{f}_{n+1}\|^2
     \, .
\end{align*}
Rearranging terms leads to
\begin{align*}
&\|\ve{u}_{n+1}\|^2-\|\ve{u}_{n}\|^2 + \tau\Big(\big\|\sqrt{\nu_{n}}\,\nabla\ve{u}_{n+1}\big\|^2-\nu_{\mathrm{min}}\|\nabla\ve{u}_{n}\|^2\Big) + \tau\big\|\sqrt{\nu_{n}}\,\nabla\ve{u}_{n+1}\big\|^2 \\   \leq 
& \left(\frac{\varepsilon\tau}{T}+\frac{2\tau\|\nabla\nu\|_{\infty}^2}{\nu_{\mathrm{min}}}\right)(\|\ve{u}_{n}\|^2 + \tau\nu_{\mathrm{min}}\|\nabla\ve{u}_{n}\|^2) + \left(2\tau+\frac{T}{\varepsilon}\right)\tau\|\ve{f}_{n+1}\|^2\, .
\end{align*}
Adding up from $n=0$ to $n=N-1$ gives
\begin{align*}
&\|\ve{u}_{N}\|^2+\tau\sum_{n=0}^{N-1}\Big(\big\|\sqrt{\nu_{n}}\,\nabla\ve{u}_{n+1} \big\|^2 -\nu_{\mathrm{min}}\|\nabla\ve{u}_{n} \|^2\Big)+\sum_{n=1}^{N}\tau\big\|\sqrt{\nu_{n-1}}\,\nabla\ve{u}_{n}\big\|^2\\  \leq \ & \|\ve{u}_0\|^2 + \Big(\frac{\varepsilon\tau}{ T} + \frac{2\tau\|\nabla\nu\|_{\infty}^2}{\nu_{\mathrm{min}}}\Big)\sum_{n=0}^{N-1}\left(\|\ve{u}_{n}\|^2 + \tau\nu_{\mathrm{min}}\|\nabla\ve{u}_{n}\|^2\right) + \left(2\tau+\frac{T}{\varepsilon} \right)\sum_{n=1}^{N}\tau\|\ve{f}_{n}\|^2\, .
\end{align*}
Finally, using 
\begin{align*}
\sum_{n=0}^{N-1}\Big(\big\|\sqrt{\nu_{n}}\,\nabla\ve{u}_{n+1} \big\|^2 -\nu_{\mathrm{min}}\|\nabla\ve{u}_{n} \|^2\Big) 
\geq \big\|\sqrt{\nu_{N-1}}\,\nabla\ve{u}_{N} \big\|^2 - \nu_{\mathrm{min}}\|\nabla\ve{u}_{0}\|^2\, ,
\end{align*}
estimate \eqref{mainResult} follows directly from the discrete Gronwall lemma.}

With the stability result just proven, we eliminate the need to know the viscosity at time step $n+1$ and \textit{still} have an unconditionally stable scheme. The price to pay is a higher regularity requirement on the viscosity field, which must now be Lipschitz continuous in space. This higher requirement may not be satisfied in the fully discrete case for some practical scenarios, but in Section~\ref{sec_nonNewtonian} we will discuss possible remedies.

\section{Fractional-step schemes}\label{sec_incremental}

This section is devoted to extending the IMEX methods from Section~\ref{sec_mono} to a class of fractional-step schemes: the incremental pressure correction method. We shall consider this prototypical scheme in its standard variant, which allows for variable viscosity \cite{Deteix2019}.

\subsection{IMEX stress-divergence formulation}
For the Navier--Stokes equations in stress-divergence form, consider the following first-order IMEX fractional-step scheme: find the solution $(\ve{u}_{n+1},\hat{\ve{u}}_{n+1},p_{n+1})$ of the system
\begin{align*}
\frac{1}{\tau}(\ve{u}_{n+1}-\hat{\ve{u}}_{n})-\nabla\cdot(\nu_{n+1}\nabla\ve{u}_{n+1}) + \mathbf{c}(\ve{u}_n,\ve{u}_{n+1})
&= \nabla\cdot(\sqrt{\nu_{n+1}}\sqrt{\nu_{n}}\,\nabla^{\top}\ve{u}_{n}) - \nabla p_{n} + \ve{f}_{n+1}\,, \\
\frac{1}{\tau}(\hat{\ve{u}}_{n+1}-\ve{u}_{n+1})+\nabla(p_{n+1} - p_n)  &= \ve{0}\,, 
\\
\nabla\cdot\hat{\ve{u}}_{n+1} & = 0\, , 
\end{align*}
with boundary conditions $\ve{u}_{n+1} = \ve{0}$ and $\ve{n}\cdot\hat{\ve{u}}_{n+1} = 0$ on $\Gamma$.
Since $\ve{u}$ is no longer divergence-free, the skew-symmetrised form $\mathbf{c}$ of the convective term \eqref{skewSymmetrisation} is used. In addition, since the initial condition $\ve{u}_0$ is divergence-free, we assume $p_0=0$ for the rest of this section.

We can use simple algebraic manipulation to eliminate the end-of-step velocity $\hat{\ve{u}}$ from the system, reformulating the $(\hat{\ve{u}},p)$-subproblem as a pressure Poisson equation. The resulting weak formulation is: find $(\ve{u}_{n+1},p_{n+1})\in [H^1_0(\Omega)]^d\times [L^2_0(\Omega)\cap H^1(\Omega)]$ such that
\begin{equation}
\begin{split}
\frac{1}{\tau}(\ve{u}_{n+1}-\ve{u}_{n},\ve{w})+(\nu_{n+1}\nabla\ve{u}_{n+1},\nabla\ve{w}) +((\nabla\ve{u}_{n+1})\ve{u}_{n},\ve{w}) + \frac{1}{2}((\nabla\cdot\ve{u}_n)\ve{u}_{n+1},\ve{w}) &=
  \\  (2p_{n}-p_{n-1},\nabla\cdot\ve{w})-\left(\sqrt{\nu_n}\,\nabla^{\top}\ve{u}_{n},\sqrt{\nu_{n+1}}\,\nabla\ve{w}\right) + (\ve{f}_{n+1},\ve{w})\,, &  \\
\tau(\nabla(p_{n+1}-p_{n}),\nabla q) + \left(q,\nabla\cdot\ve{u}_{n+1}\right) &= 0\,,
\label{fractionalSD}
\end{split}
\end{equation}
for all $(\ve{w},q)\in [H^1_0(\Omega)]^d\times [L^2_0(\Omega)\cap H^1(\Omega)]$. This split-step scheme is also unconditionally stable, as the following result shows. 

\begin{theorem}[Stability of an IMEX fractional-step scheme in SD form] 
Let the time step be given by $\tau=T/N$ with $N\ge 2$. 
Then, the fractional-step method \eqref{fractionalSD} satisfies the following stability estimate
\begin{equation}
\begin{split}
    \|\ve{u}_{N}\|^2 + \tau\big\| \sqrt{\nu_{N}}\,\nabla\ve{u}_{N}\big\|^2 + \tau^2\|\nabla p_N\|^2 + \frac{\tau}{T}\sum_{n=1}^{N}\Big(\tau\big\| \sqrt{\nu_{n}}\,\nabla\ve{u}_{n}\big\|^2+\tau^2\|\nabla p_n\|^2\Big) \\ \leq M\Bigg[\|\ve{u}_{0}\|^2 + \tau\|\sqrt{\nu_0}\,\nabla\ve{u}_0 \|^2 + T\sum_{n=1}^{N}\tau\|\ve{f}_{n}\|^2\Bigg]\,, 
\label{stabilityMonoSqrtIncremental}
\end{split}
\end{equation}
where
\begin{align*}
M &= \mathrm{exp}\Bigg(\frac{1}{1-\tau/T}\Bigg) = \mathrm{exp}\Bigg(\frac{N}{N-1}\Bigg) \leq \mathrm{e}^2\, .
\end{align*}
\end{theorem}
\proof{Setting $\ve{w}=2\tau\ve{u}_{n+1}$ in \eqref{fractionalSD} yields
\begin{equation}
\begin{split}
     &\|\ve{u}_{n+1}\|^2 + \|\delta\ve{u}_{n+1}\|^2 - \|\ve{u}_{n}\|^2  + 2\tau\big\|\sqrt{\nu_{n+1}}\,\nabla\ve{u}_{n+1}\big\|^2 \\ &=  2\tau\left[\left(2p_n-p_{n-1},\nabla\cdot\ve{u}_{n+1}\right)-\left(\sqrt{\nu_{n}}\,\nabla^{\top}\ve{u}_{n},\sqrt{\nu_{n+1}}\,\nabla\ve{u}_{n+1}\right) + (\ve{f}_{n+1},\ve{u}_{n+1})\right] \\ 
    &\leq 2\tau\left[\left(2p_n-p_{n-1},\nabla\cdot\ve{u}_{n+1}\right)+\big\|\sqrt{\nu_{n}}\,\nabla\ve{u}_{n}\big\|^2+\big\|\sqrt{\nu_{n+1}}\,\nabla\ve{u}_{n+1}\big\|^2 + (\ve{f}_{n+1},\ve{u}_{n+1})\right].
\label{auxResultFractional}
\end{split}
\end{equation}

Compared to the monolithic cases presented in the last section, we have an additional pressure term to estimate. Denoting $\delta p_{n+1}:=p_{n+1}- p_{n}$, the second equation in \eqref{fractionalSD} implies 
\begin{equation}
\tau(\nabla\delta p_{n+1} - \nabla\delta p_{n},\nabla q) = -(\nabla\cdot\delta\ve{u}_{n+1},q) = \left(\delta\ve{u}_{n+1},\nabla q\right)\,.
\label{continuityDelta}
\end{equation}
Now, setting $q = 2\tau(2p_{n}-p_{n-1})$ in \eqref{fractionalSD} gives
\begin{align*}
2\tau(2p_n-p_{n-1},\nabla\cdot\ve{u}_{n+1}) &= -2\tau^2(\nabla(2p_n-p_{n-1}),\nabla(p_{n+1}-p_n))\\
&\equiv -2\tau^2\big(\nabla[p_{n+1}-(\delta p_{n+1}-\delta p_{n})],\nabla(p_{n+1}-p_n)\big)\\
&= -2\tau^2(\nabla p_{n+1},\nabla p_{n+1}-\nabla p_n) + 2\tau^2(\nabla\delta p_{n+1}-\nabla\delta p_n,\nabla\delta p_{n+1}).
\end{align*}
For the first term on the right-hand side above we have
\begin{align*}
    -2(\nabla p_{n+1},\nabla p_{n+1}-\nabla p_n) = \|\nabla p_{n}\|^2-\|\nabla p_{n+1}\|^2-\|\nabla\delta p_{n+1}\|^2\, ,
\end{align*}
and for the second one we set $q=\tau\delta p_{n+1}$ in \eqref{continuityDelta}, so that
\begin{align*}
     2\tau^2(\nabla\delta p_{n+1}-\nabla\delta p_n,\nabla\delta p_{n+1}) = 2\tau(\delta\ve{u}_{n+1},\nabla\delta p_{n+1}) \leq \|\delta\ve{u}_{n+1}\|^2 + \tau^2\|\nabla\delta p_{n+1}\|^2\,. 
\end{align*}
Collecting these last results, we estimate the pressure term in \eqref{auxResultFractional} as
\begin{align}
2\tau(2p_n-p_{n-1},\nabla\cdot\ve{u}_{n+1}) \leq \tau^2\|\nabla p_{n}\|^2-\tau^2\|\nabla p_{n+1}\|^2 + \|\delta\ve{u}_{n+1}\|^2\,.\label{Two-star}
\end{align}
Moreover,
\begin{align}
    2\tau(\ve{f}_{n+1},\ve{u}_{n+1}) \leq \frac{\tau}{T}\|\ve{u}_{n+1}\|^2+T\tau\|\ve{f}_{n+1}\|^2\,.\label{Three-star}
\end{align}
Replacing \eqref{Two-star} and \eqref{Three-star} in \eqref{auxResultFractional} we get
\begin{align*}
     &\left(\|\ve{u}_{n+1}\|^2 + \big\|\sqrt{\nu_{n+1}}\,\nabla\ve{u}_{n+1}\big\|^2+\tau^2\|\nabla p_{n+1}\|^2\right)  - \left(\|\ve{u}_{n}\|^2 + \big\|\sqrt{\nu_{n}}\,\nabla\ve{u}_{n}\big\|^2+\tau^2\|\nabla p_{n}\|^2\right) \\  
    &\leq \frac{\tau}{T}\|\ve{u}_{n+1}\|^2+T\tau\|\ve{f}_{n+1}\|^2\, ,
\end{align*}
and adding from $n=1$ to $n=N-1$ leads to
\begin{align*}
     &\|\ve{u}_{N}\|^2 + \tau\big\|\sqrt{\nu_{N}}\,\nabla\ve{u}_{N}\big\|^2+\tau^2\|\nabla p_{N}\|^2  \\ &\leq  \|\ve{u}_{0}\|^2 + \tau\big\|\sqrt{\nu_{0}}\,\nabla\ve{u}_{0}\big\|^2+\tau^2\|\nabla p_{0}\|^2 +\frac{\tau}{T}\sum_{n=1}^{N}\|\ve{u}_{n}\|^2+T\sum_{n=1}^{N}\tau\|\ve{f}_{n}\|^2\,.
\end{align*}
Finally, adding
\begin{align*}
    \frac{\tau}{T}\sum_{n=1}^{N}\Big(\tau\big\|\sqrt{\nu_{n}}\,\nabla\ve{u}_{n}\big\|^2+\tau^2\|\nabla p_{n}\|^2\Big)
\end{align*}
to both sides of the last inequality, and noticing that $\tau/T<1$, we conclude by using the Gronwall inequality from Lemma~\ref{Lem:Gronwall-conditional}.}

As it was the case with Theorem~\ref{Th:Monosqrt}, the stability result just shown is in a weak norm, and also requires the knowledge of $\nu_{n+1}$. To remedy this, we will present below an IMEX fractional-step scheme based on the generalised Laplacian rewriting of the diffusive term.

\subsection{IMEX generalised Laplacian formulation}
After elimination of the end-of-step velocity, the weak problem for the GL version of the IMEX incremental projection scheme is: find
$(\ve{u}_{n+1},p_{n+1})\in [H^1_0(\Omega)]^d\times [L^2_0(\Omega)\cap H^1(\Omega)]$ such that
\begin{equation}
\begin{split}
\frac{1}{\tau}(\ve{u}_{n+1}-\ve{u}_{n},\ve{w})+\left(\nu_{n}\nabla\ve{u}_{n+1},\nabla\ve{w}\right) +((\nabla\ve{u}_{n+1})\ve{u}_{n},\ve{w}) + \frac{1}{2}((\nabla\cdot\ve{u}_n)\ve{u}_{n+1},\ve{w}) &=
  \\  \left(2p_{n}-p_{n-1},\nabla\cdot\ve{w}\right)+\left(\nabla^{\top}\ve{u}_{n}\nabla\nu_n,\ve{w}\right) + \left( \ve{w},\ve{f}_{n+1} \right)\,,&\\
\tau(\nabla(p_{n+1}-p_{n}),\nabla q) + \left(\nabla\cdot\ve{u}_{n+1},q\right) &= 0\,,
\label{fractionalGL}
\end{split}
\end{equation}
for all $(\ve{w},q)\in [H^1_0(\Omega)]^d\times [L^2_0(\Omega)\cap H^1(\Omega)]$. The following result states the stability for this scheme. 

\begin{theorem}[Stability of the IMEX fractional-step scheme in GL form]
Let \\ us assume that $\nabla\nu\in [L^\infty(\bar{Q}]^d$. Then, the fractional-step IMEX-GL scheme \eqref{fractionalGL} yields the stability estimate
\begin{equation}
\begin{split}
      \|\ve{u}_{N}\|^2   +\gamma\tau\big\|\sqrt{\nu_{N-1}}\,\nabla \ve{u}_{N}\big\|^2 + \tau^2\|\nabla p_{N}\|^2 + \sum_{n=1}^{N}\tau\Big[(2-\gamma)\big\|\sqrt{\nu_{n-1}}\,\nabla\ve{u}_{n}\big\|^2 +\alpha\tau^2\|\nabla p_{n}\|^2\Big] \\ \leq  \mathrm{exp}\left(\frac{\alpha T}{1-\alpha\tau}\right)\left[\|\ve{u}_0 \|^2 + \tau\nu_{\mathrm{min}}\|\nabla\ve{u}_0\|^2   + \varepsilon T\sum_{n=1}^{N}\tau\| \ve{f}_n\|^2\right]\,,\label{conditionalStability}
\end{split}
\end{equation}
 for any positive time-step size $\tau$ satisfying 
 \begin{equation}
    \tau < \left(\frac{\|\nabla\nu\|_{\infty}^2}{\gamma\nu_{\mathrm{min}}}+\frac{1}{\varepsilon T}\right)^{-1} =: \alpha^{-1} \, ,
    \label{conditionDeltaT}
 \end{equation}  
 with $\gamma \in (0,2)$ and $\varepsilon > 0$, both chosen arbitrarily.
\end{theorem}
\proof{Compared to the stress-divergence version, the only term estimated differently is
\begin{align*}
    2\tau(\nabla^{\top}\ve{u}_n\nabla\nu_n + \ve{f}_{n+1},\ve{u}_{n+1}) &\leq 2\tau\|\nabla\nu\|_{\infty}\|\nabla\ve{u}_n\|\,\|\ve{u}_{n+1}\| + 2\tau\|\ve{f}_{n+1}\|\,\|\ve{u}_{n+1}\|\\
    &\leq \gamma\tau\nu_{\mathrm{min}}\|\nabla\ve{u}_n\|^2+\Big(\frac{\tau\|\nabla\nu\|_{\infty}^2}{\gamma\nu_{\mathrm{min}}}+\frac{\tau}{\varepsilon T}\Big)\|\ve{u}_{n+1}\|^2 + \tau\varepsilon T\|\ve{f}_{n+1}\|^2 ,  
\end{align*}
in which $0<\gamma < 2$ and $\varepsilon > 0$. Hence,
\begin{align*}
\|\ve{u}_{n+1}\|^2 + \tau^2\|\nabla p_{n+1}\|^2  + 2\tau\big\|\sqrt{\nu_{n}}\,\nabla\ve{u}_{n+1}\big\|^2 \\ \leq \|\ve{u}_{n}\|^2 + \tau^2\|\nabla p_{n}\|^2 + \gamma\tau\nu_{\mathrm{min}}\|\nabla\ve{u}_{n}\|^2 +\Big(\frac{\tau\|\nabla\nu\|_{\infty}^2}{\gamma\nu_{\mathrm{min}}}+\frac{\tau}{\varepsilon T}\Big)\|\ve{u}_{n+1}\|^2 +  \tau\varepsilon T\|\ve{f}_{n+1}\|^2 .
\end{align*}
Adding up from $n=0$ to $n=N-1$ and using that $p_0=0$ yields 
\begin{equation}
\begin{split}
& \|\ve{u}_{N}\|^2+\gamma\tau\big\|\sqrt{\nu_{N-1}}\,\nabla\ve{u}_{N}\big\|^2 + \tau^2\|\nabla p_{N}\|^2 + (2-\gamma)\tau\sum_{n=1}^{N}\big\|\sqrt{\nu_{n-1}}\,\nabla\ve{u}_{n}\big\|^2\\  \leq \ &
 \|\ve{u}_{0}\|^2+\gamma\tau\nu_{\mathrm{min}}\|\nabla\ve{u}_{0} \|^2 + \tau\varepsilon T\sum_{n=1}^{N}\|\ve{f}_{n}\|^2+\Big(\frac{\tau\|\nabla\nu\|_{\infty}^2}{\gamma\nu_{\mathrm{min}}}  +\frac{\tau}{\varepsilon T}\Big)\sum_{n=1}^{N}\|\ve{u}_{n}\|^2 \, .
 \label{inequalityNplus1}
 \end{split}
\end{equation}
Finally, by adding 
\begin{align*}
   \Big(\frac{\tau\|\nabla\nu\|_{\infty}^2}{\gamma\nu_{\mathrm{min}}}+\frac{\tau}{\varepsilon T}\Big)\sum_{n=1}^{N}\tau^2\|\nabla p_{n}\|^2
\end{align*}
to both sides of \eqref{inequalityNplus1} and choosing $\tau$ according to \eqref{conditionDeltaT}, we get the stability estimate \eqref{conditionalStability} by applying Lemma~\ref{Lem:Gronwall-conditional}.}

\begin{remark}
The last stability result requires a CFL condition on the time-step size $\tau$. It is worth noticing that this condition includes only the viscosity field (i.e., the problem data). In addition, since the constants $\gamma$ and $1/\varepsilon$ in \eqref{conditionDeltaT} can be chosen arbitrarily close to $2$ and $0$, respectively, this restriction can be simplified to
\begin{equation}
    \tau < \frac{2\nu_{\mathrm{min}}}{\|\nabla\nu \|_{\infty}^2}\,,
    \label{condition}
 \end{equation}
which can also be obtained by letting $T\rightarrow\infty$, that is, letting the problem tend to a steady-state. Interestingly, \eqref{condition} is related to the condition imposed on the reaction term in \cite{Anaya2021} to analyse a mixed finite element method for the linear, stationary Oseen equation.
 \end{remark}

\section{Summary of theoretical results}\label{sec_summary}
Let us briefly collect the results of our analysis. We have proved the stability of several IMEX discretisations, and the corresponding results are summarised on Table \ref{table_stability}. For four discretisations we have unconditional stability (albeit in different norms, some stronger than others), and for another one we have proved stability under a time-step restriction \eqref{condition}. The only method where stability has not been proven is the SD formulation with fully explicit viscosity, that is, \eqref{SDnaive}. Therefore, the main purpose of our numerical examples will be to verify whether our analysis is sharp with respect to condition \eqref{condition} and to the (lack of) stability of the ``naive'' IMEX-SD discretisation.
\begin{table}[ht!]
 \centering
 \caption{Theoretical results on the stability of different IMEX treatments of the viscous term.}
 {\begin{tabular}{|c|l|l|c|}
    \hline
   \textbf{Type}   & \textbf{Left-hand side term} & \textbf{Right-hand side term} & \textbf{Stability}\\
   \hline
   Monolithic   & $-\nabla\cdot(2\nu_n\nabla^{\mathrm{s}}\ve{u}_{n+1})$ & & unconditional\\
   
   Monolithic  & $-\nabla\cdot\big(\nu_n\nabla\ve{u}_{n+1})$ & $\nabla\cdot(\nu_n\nabla^{\top}\ve{u}_{n}\big)$ & not proved\\
   
   Monolithic  & $-\nabla\cdot(\nu_{n+1}\nabla\ve{u}_{n+1})$ & $\nabla\cdot\big(\sqrt{\nu_{n+1}}\sqrt{\nu_n}\, \nabla^{\top}\ve{u}_{n}\big)$ & unconditional \\
   
   Monolithic  &  $-\nabla\cdot(\nu_n\nabla\ve{u}_{n+1})$ & $\nabla^{\top}\ve{u}_{n}\nabla\nu_n$ & unconditional \\

    Fractional  & $-\nabla\cdot(\nu_{n+1}\nabla\ve{u}_{n+1})$ & $\nabla\cdot\big(\sqrt{\nu_{n+1}}\sqrt{\nu_n}\, \nabla^{\top}\ve{u}_{n}\big)$ & unconditional \\

      Fractional    &  $-\nabla\cdot(\nu_n\nabla\ve{u}_{n+1})$ & $\nabla^{\top}\ve{u}_{n}\nabla\nu_n$ & $\tau < \tau_{\mathrm{max}}(\nu)$ \\
    \hline
 \end{tabular}}
  \label{table_stability} 
\end{table}

\section{On the fully discrete problem}\label{sec_spatial}
Although we restrict our stability analysis to the semi-discrete setting, some brief comments should be drawn on the fully discrete problem. We will restrict our discussion to conforming finite element spaces with continuous pressures, since that choice does not modify the matrix structure of the semi-discrete problems.

\subsection{Matrix structure of the IMEX methods}

One of the most appealing computational features of the IMEX schemes presented herein is the resulting matrix structure. Consider a generic spatial discretisation applied to $\ve{u}_n$ and $p_n$, denoting by $\ve{U}_{n}$ and $\ve{P}_{n}$ the respective vectors of degrees of freedom. Since we are interested in problems where the viscosity may not be known, and also interested in IMEX schemes leading to solving only linear systems at each time step, we  shall only discuss the schemes that do not require $\nu_{n+1}$. 

In the monolithic case (e.g., the scheme \eqref{strongMomentumMono}), the fully-discrete linear system to be solved at each time step has the form
\begin{align*}
    \begin{bmatrix}
        \mathbf{A}_{n} & \mathbf{B}^{\top}\, \\
        \mathbf{B} & \mathbf{0}
    \end{bmatrix}
    \begin{bmatrix}
        \ve{U}_{n+1} \\ \ve{P}_{n+1}
    \end{bmatrix} = 
    \begin{bmatrix}
        \ve{F}_{n+1} \\ \ve{0}
    \end{bmatrix},
\end{align*}
where $\mathbf{B}$ is (minus) the divergence matrix and $\ve{F}_{n+1}$ is the right-hand side vector depending on $\ve{f}_{n+1}$, $\ve{U}_{n}$, and $\nu_{n}$; the velocity-velocity matrix $\mathbf{A}_n$ will consist of $d$ identical blocks:
\begin{align}
\mathbf{A}_{n} = 
    \begin{bmatrix}
        \mathbf{K}_{n} & \ve{0}   \\
        \ve{0} & \mathbf{K}_{n}
    \end{bmatrix} 
    \  \text{for} \  d=2\, , \ \, \text{or} \ 
\mathbf{A}_{n} = 
    \begin{bmatrix}
        \mathbf{K}_{n} & \ve{0}  & \ve{0}  \\
         \ve{0} & \mathbf{K}_{n} & \ve{0} \\
         \ve{0} & \ve{0} & \mathbf{K}_{n}
    \end{bmatrix} 
    \  \text{for} \  d=3\, ,
    \label{blockDiagonal}
\end{align}
where $\mathbf{K}_{n} := \mathbf{M} + \mathbf{D}(\nu_n) + \mathbf{C}(\ve{U}_{n})$, with $\mathbf{M}$, $\mathbf{D}$ and $\mathbf{C}$ denoting, respectively, the standard mass, diffusion, and convection matrices from a scalar transport equation. The same structure will be observed for higher-order temporal discretisations, with $\nu_n$ and $\ve{U}_{n}$ replaced by their corresponding extrapolations. This block structure reduces the costs of assembling and solving the algebraic system \cite{John2016}, being especially advantageous for fractional-step schemes, as it allows us to split the velocity update into $d$ scalar problems:
\begin{align*}
    \mathbf{K}_n\, \ve{U}^i_{n+1} = \ve{F}^i(\ve{P}_{n},\ve{P}_{n-1},\ve{U}_{n},\nu_n)\, ,\ \ i=1,...\, ,d\, ,
\end{align*}
where $\ve{U}^i_{n+1}$ contains the degrees of freedom of the $i$-th velocity component. After computing all the components, the pressure is updated by solving another scalar problem.

\begin{remark}
One may also wish to treat convection explicitly, to eliminate from $\mathbf{K}_n$ the non-symmetric matrix $\mathbf{C}$, so that more efficient solvers can be used. The price to pay for this choice is, in most cases, a CFL condition, as it was the case in recent works \cite{Burman2022,Burman2023}.
\end{remark}

\subsection{Generalised Newtonian fluids}\label{sec_nonNewtonian}

So far in our analysis we have considered that the viscosity $\nu$ is fixed and known. Nevertheless, in several interesting applications (such as hemodynamic flows and nanofluids, see e.g. \cite{ABUGATTAS2020529,FARIAS2023,GONZALEZ2021101635,AGUIRRE2022104400}) the fluid is best described using a generalised Newtonian constitutive law. In such a case, the viscosity is given by
\begin{align*}
    \nu(\ve{x},t) = \eta(|\nabla^{\mathrm{s}}\ve{u}(\ve{x},t)|)\, ,
\end{align*}
where $|\cdot|$ denotes the Euclidean norm, and $\eta$ is a continuous function. A popular example is the Carreau model \cite{Galdi2008}
\begin{align}
    \eta(|\nabla^{\mathrm{s}}\ve{u}|) = \nu_{\infty} + (\nu_{0}-\nu_{\infty})(1+2\lambda^2|\nabla^{\mathrm{s}}\ve{u}|^2)^{-m}\, , 
    \label{Carreau}
\end{align}
where $\nu_0$, $\nu_\infty$ and $\lambda$ are positive, experimental rheological parameters. For the common case of shear-thinning fluids ($\nu_{\infty} > \nu_{0}$), we get $\nu_{0} \leq \nu \leq \nu_{\infty}$.

The viscosity just described, although covered by our analysis, may cause issues in the fully discrete case, where discontinuities in the velocity gradient lead to discontinuous viscosity. Several approaches can be used to overcome this. For example, local averaging can be used to produce a continuous viscosity field. Alternatively, computing the $L^2(\Omega)$-orthogonal projection onto a space of continuous finite element functions would also fit our hypotheses. This last alternative has been used, with satisfactory numerical results, for schemes related to the ones presented herein \cite{Schussnig2021JCP,Pacheco2021CMAME}, and it is also the approach we follow in our numerical experiments. The mathematical analysis of the resulting schemes has not been carried out so far, and it is out of the scope of the present paper.


\section{Numerical examples}\label{sec_examples}

In this section we present three series of numerical experiments testing the difference in stability properties for the IMEX schemes. In accordance to the discussion in the last section, we have considered conforming Taylor-Hood elements in quadrilateral meshes (biquadratic velocities and continuous bilinear pressures). For problems with a velocity-dependent viscosity we have projected the viscosity onto the continuous bilinear finite element space. Throughout this section, whenever an IMEX-SD formulation is mentioned, we will be referring to one for which we have proved unconditional stability -- unless where otherwise stated.


\subsection{Temporal convergence}

\begin{figure}[ht!]
\centering
\includegraphics[trim = 55 0 77 0,clip, width = 1.00\textwidth]{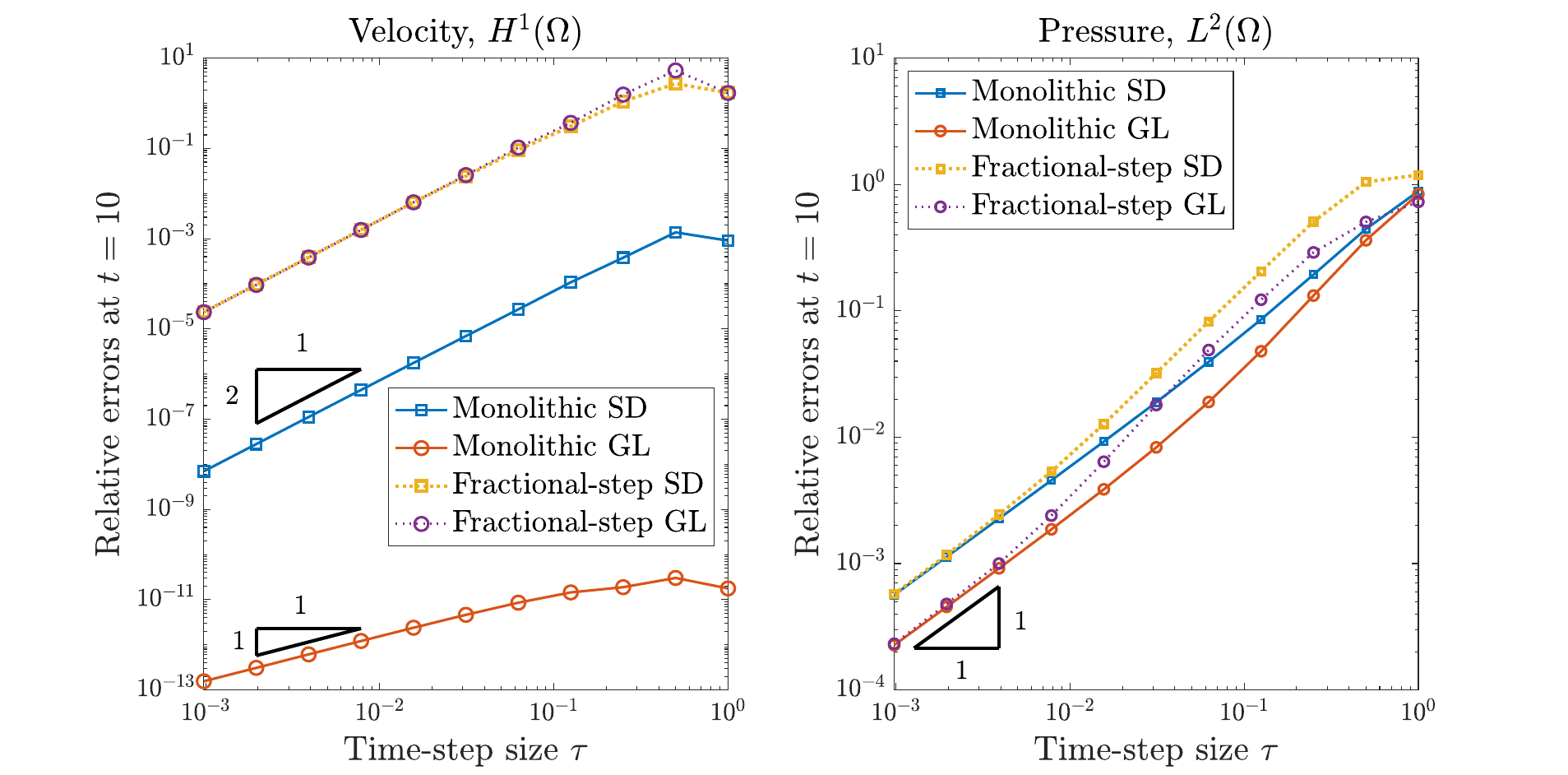}
\caption{Temporal convergence study for the IMEX schemes. The quadratic velocity convergence of the fractional-step methods is probably an initial superconvergence that would break down with further temporal refinement.}
\label{temporalConvergence}
\end{figure}

Our first test aims at assessing the temporal accuracy/convergence of the IMEX schemes presented herein. For this purpose, we use a problem with known exact solution. The viscosity field is given by
\begin{equation}
    \nu(x,y,t) = xyf(t) + g(t)\, ,
\end{equation}
where $f$ and $g$ are any two functions that guarantee $\nu>0$ in $\Omega$ for all $t\in [0,T]$. With this choice, the exact solution of the Navier--Stokes equation is given by 
\begin{align*}
    \ve{u} = 2f(t)
    \begin{pmatrix}
    y\\x
    \end{pmatrix} \quad \text{and} \quad p = (C-2xy)f'(t)\, ,
\end{align*}
 with $\ve{f}\equiv\ve{0}$ and the boundary conditions chosen  accordingly; the constant $C$ is set so that $p\in L^2_0(\Omega)$ for all $t$. We have chosen this solution since, in space, it belongs to the finite element space, so the entire error stems from time discretisation. We set $\Omega = (0,1)^2$, $T=10$, $f(t)=\sin ^2 t$ and $g(t) \equiv 0.001$. The square domain is discretised uniformly with $4\times4$ square elements, which completely resolves (at each time instant) the exact solution in space. The time-step size is halved at each refinement level, starting from $\tau=1$. In Figure~\ref{temporalConvergence} we depict the errors in velocity and pressure at $t=10$. For all schemes considered, the convergence is at least linear in both fields, with the monolithic GL being the most accurate for velocity (although the remaining schemes show a higher order of convergence), and a comparable performance of all methods for pressure.
 
\begin{remark}
It is interesting to notice that all methods present a very stable behaviour even for the larger time step sizes. This, despite the fact that the viscosity field considered in this example is such that $\nu_{\mathrm{min}}/\|\nabla \nu \|^2_{\infty} = 5\times 10^{-4}$. According to our analysis, this would require $\tau < 0.001$ to guarantee the stability of the fractional-step GL method. That, however, was not observed in this study.
\end{remark}

\subsection{Stability test}\label{sec_cavity}
To put the stability of the fractional-step IMEX schemes to a more challenging test, we tackle a problem with less smooth solution. The lid-driven cavity flow is one of the most popular benchmarks for generalised Newtonian flow solvers, and a common setup involves a power-law fluid:
\begin{align}
    \eta(|\nabla^{\mathrm{s}}\ve{u}|) = \kappa\big|2^{\frac 12}\nabla^{\mathrm{s}}\ve{u}\big|^{n-1}  \, ,
    \label{powerLaw}
\end{align}
where $\kappa$, $n$ and $\nu_{\mathrm{min}}$ are positive constants. For the square cavity $\Omega = (0,1)^2$, the boundary conditions are a horizontal velocity $U(t)$ on the lid ($y=1$), and no slip on the bottom and lateral walls. We consider a reference shear-thickening fluid with $\kappa=0.01$ and $n=1.5$ \cite{Neofytou2005}. To avoid a zero viscosity, we replace \eqref{powerLaw} by 
\begin{align*}
    \eta(|\nabla^{\mathrm{s}}\ve{u}|) = \mathrm{max}\left\lbrace 10^{-5},\, \kappa\big|2^{\frac 12}\nabla^{\mathrm{s}}\ve{u}\big|^{n-1}  \right\rbrace .
\end{align*}
The lid velocity is prescribed according to
\begin{align}
 U(t) = \frac{1}{2}\left\lbrace 1+\text{sign}(t-t^{\star})+\left[1-\text{sign}(|2t-t^{\star}|-t^{\star})\right]\sin^2\Big(\frac{\pi t}{2t^{\star}}\Big)\right\rbrace,
 \label{ramp}
\end{align}
which ramps up smoothly from $U=0$ (at $t=0$) to $U=1.0$ (for $t\geq t^{\star}=2$). We run the simulation up to $T=30$, aiming for a steady-state flow that can be compared to reference stationary results \cite{Neofytou2005,Li2014}. Two uniform meshes are considered: one with $200\times 200$ elements, and one twice as fine. To test the stability of the different schemes, we take the time-step size rather large, namely $\tau = 0.25$. The aim of this test is twofold: firstly, we want to verify whether the time-step restriction \eqref{condition} found for the fractional-step GL scheme is sharp -- and we do that by tackling a problem where $\|\nabla\nu\|_{\infty}/\nu_{\mathrm{min}}\rightarrow \infty$. Second, we aim to illustrate the poor stability of IMEX-SD schemes when the coupling term $\nabla\cdot(\nu\nabla^{\top}\ve{u})$ is made fully explicit (and, since the viscosity $\nu_{n+1}$ is an unknown of the problem, we can only consider the scheme \eqref{SDnaive} as IMEX-SD). To do that, three schemes are considered:
\begin{enumerate}
   \item Fractional-step IMEX, generalised Laplacian form \eqref{fractionalGL}.
    \item Fractional-step IMEX, stress-divergence form using only $\nu_{n}$ (not $\nu_{n+1}$).
    \item Monolithic, fully coupled stress-divergence formulation \eqref{fullyCoupled}.
\end{enumerate}
The third variant is meant as a control case, since it is the method for which we have the strongest stability estimate. 

\begin{figure}[ht!]
\centering
\includegraphics[trim = 50 7 50 10,clip, width = .99\textwidth]{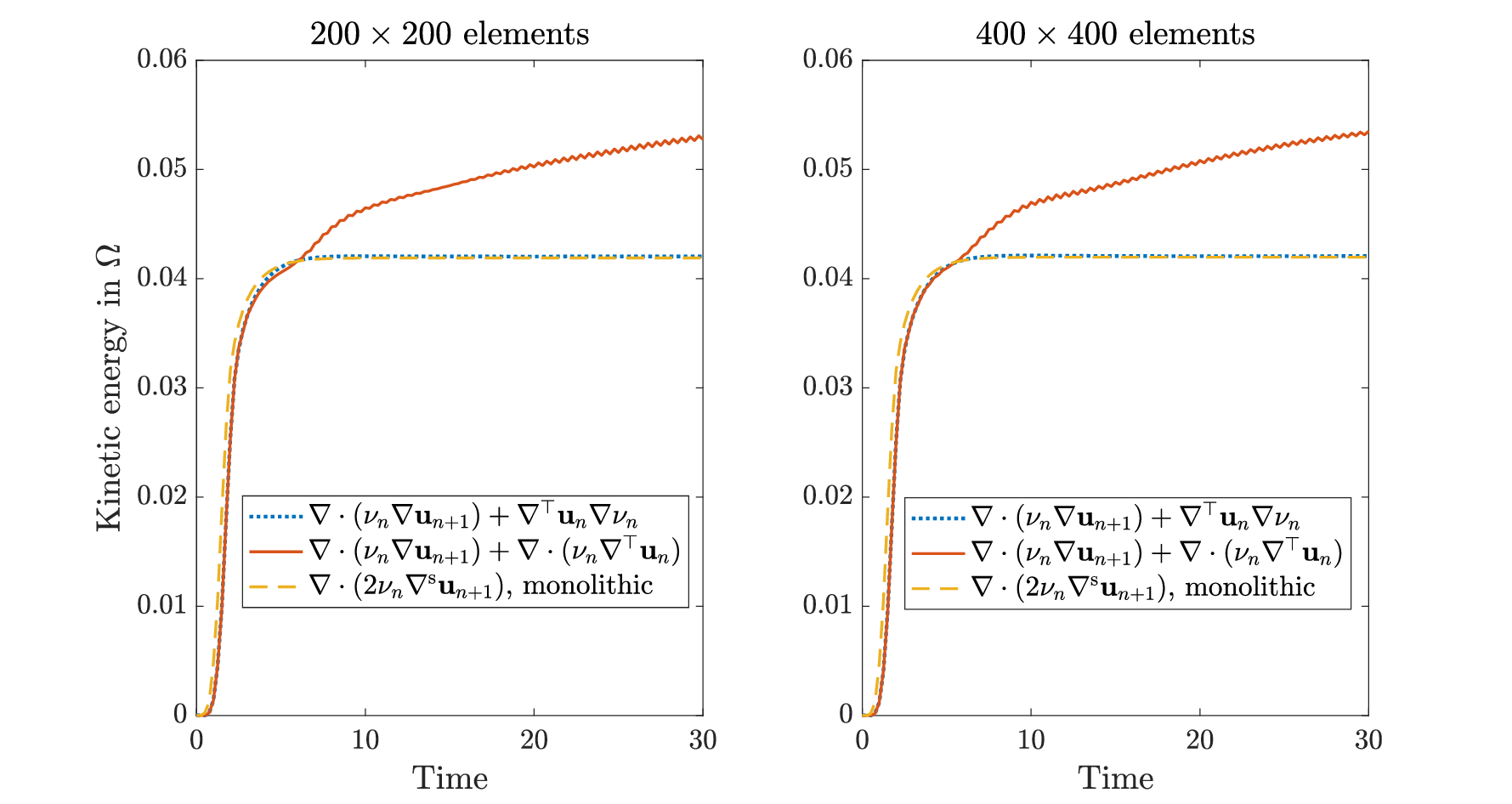}
\caption{Lid-driven cavity flow with power-law fluid: temporal evolution of the kinetic energy, confirming that a stress-divergence scheme with fully explicit viscosity can be unstable.}
\label{stabilityTest}
\end{figure}

\begin{figure}[ht!]
\centering
\includegraphics[trim = 10 50 20 50,clip, width = .8\textwidth]{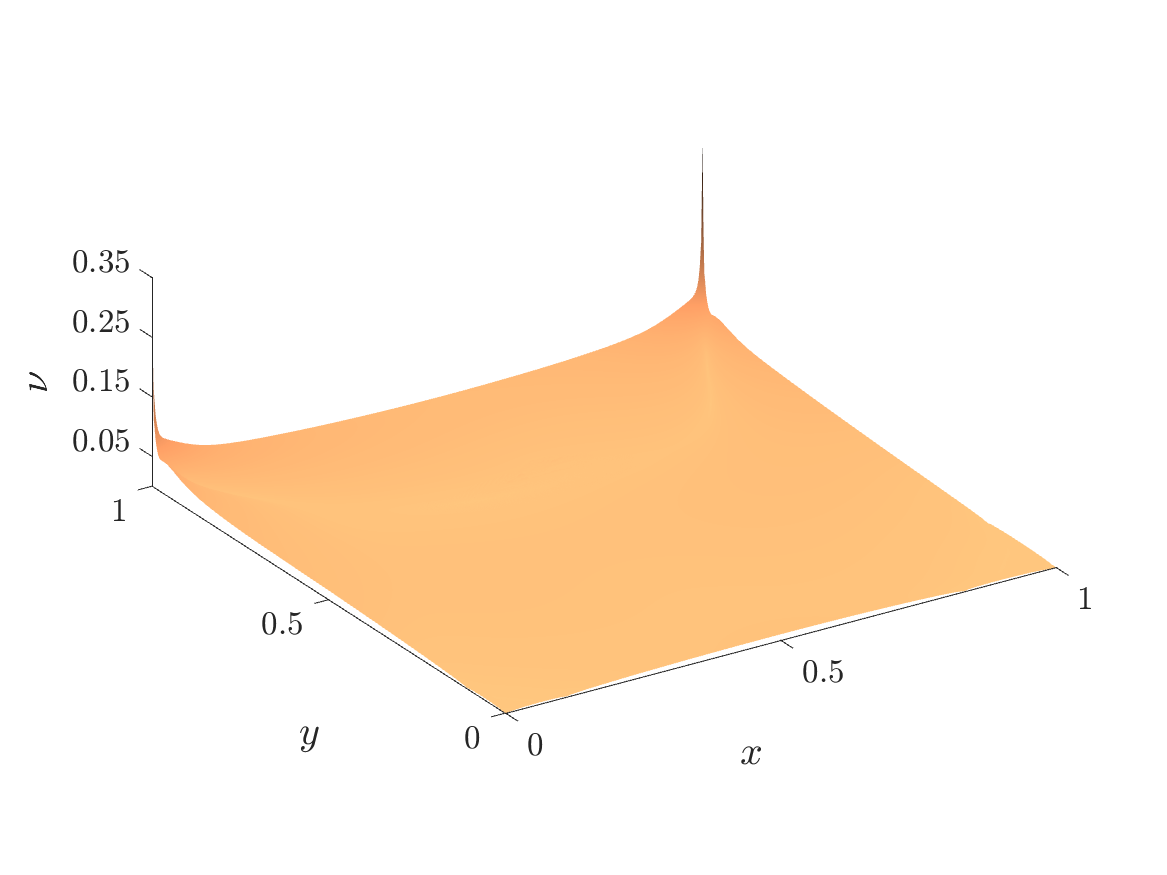}
\caption{Lid-driven cavity flow with power-law fluid: steady-state viscosity field obtained using the IMEX-GL method.}
\label{horns}
\end{figure}

To verify the accuracy of our solution, we use the location of the primary vortex as benchmark. Table \ref{table_lidVortex} shows the comparison between our results (with the coarser mesh) and reference ones \cite{Neofytou2005,Li2014}. Among the three methods we used, the two with proven stability are in good agreement with the references; on the other hand, for the IMEX-SD \eqref{SDnaive} (with fully explicit viscosity) no primary vortex position could be clearly identified. Figure \ref{stabilityTest} shows how the kinetic energy $E(t):=\|\ve{u}\|^2/2$ evolves with time: while the two stable schemes reach a steady-state level, IMEX-SD \eqref{SDnaive} exhibits a clearly unstable energy growth (notice that the spatial refinement does not remedy the instability). The two negative results for the IMEX-SD \eqref{SDnaive} are important as they hint strongly that, in general, one cannot guarantee unconditional stability of schemes such as \eqref{SDnaive}. Furthermore, Figure~\ref{horns} depicts the viscosity field at the final time, where we see the spikes near the corners (which lead to very high gradients). The time-step size used for all methods in this experiment ($\tau=0.25$) clearly violates \eqref{condition}. Yet, the fractional-step GL method again showed a stable behaviour, which hints at the possibility that \eqref{condition} is pessimistic.

\begin{table}[ht!]
 \centering
 \caption{Lid-driven cavity flow with power-law fluid: steady-state position $(x,y)$ of the primary vortex for three schemes: fully coupled SD \eqref{fullyCoupled}, fractional-step IMEX-GL \eqref{fractionalGL} and fractional-step IMEX-SD using only $\nu_n$. For the latter, the solution did not converge to a steady-state, so it was not possible to determine clear vortex locations.}
 {\begin{tabular}{|c|c|c|c|c|}
    \hline
   \textbf{Ref.~\cite{Neofytou2005}}   & \textbf{Ref.~\cite{Li2014}} & \textbf{Fully coupled} & \textbf{IMEX-GL} & \textbf{IMEX-SD, $\nu_n$}\\
   \hline
   $(0.565,0.724)$   & $(0.5628,0.7282)$ & $(0.562,0.728)$ & $(0.560,0.728)$ & undetermined\\
    \hline
 \end{tabular}}
  \label{table_lidVortex} 
\end{table}

\subsection{Idealised hemodynamic flow}
\begin{figure}[ht!]
\centering
\includegraphics[trim = 0 0 0 0,clip, width = .85\textwidth]{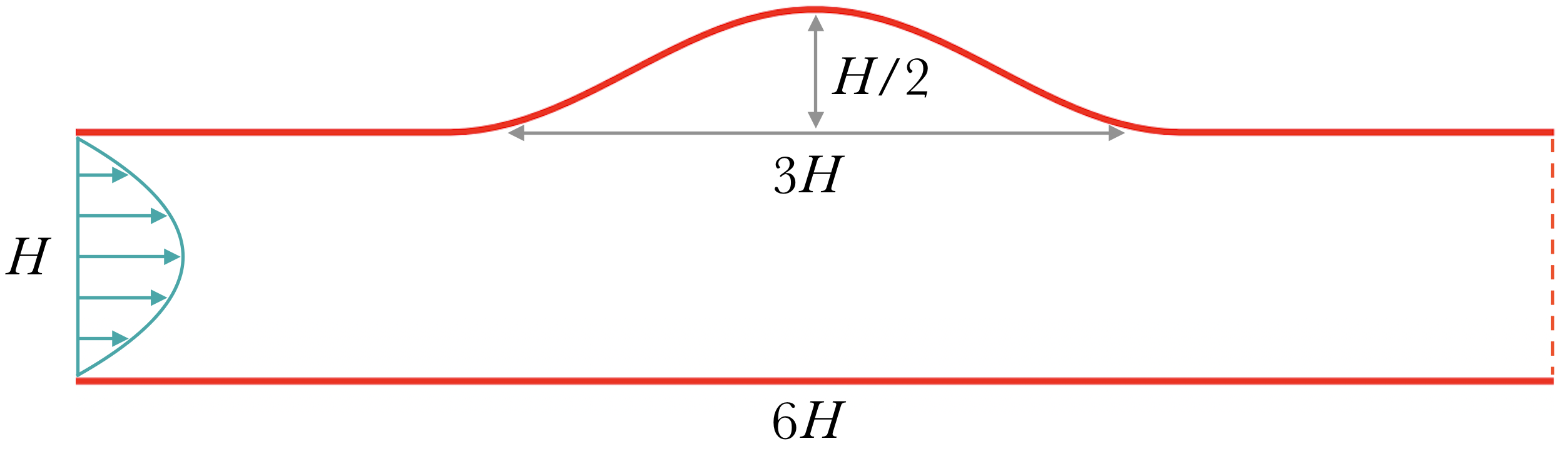}
\caption{Ccomputational domain representing an idealised ICA aneurysm.}
\label{aneurysmGeometry}
\end{figure}

Our last numerical test tackles a truly dynamic flow problem. A popular application for variable viscosity is hemodynamics, with blood usually modelled as a shear-thinning fluid. We have applied the Carreau model \eqref{Carreau} with the following parameters: $\nu_{\infty}=3.286\times 10^{-6}$ m$^2$/s, $\nu_{0}=53.33\times 10^{-6}$ m$^2$/s , $\lambda=3.313$ s and $m=0.3216$ \cite{Cho1991}. Combined with those, we consider an idealised two-dimensional aneurysm, whose geometry is illustrated in Figure \ref{aneurysmGeometry}. The upper wall is described by the smooth curve
\begin{align*}
    y(x) = H + \frac{H}{4}[1-\mathrm{sign}(|2x-6H|-3H)]\cos^2\Big(\pi\frac{x-3H}{3H}\Big)\, ,
\end{align*}
with $H = 2.5\times 10^{-3}$ m.

To simulate a physiological regime \cite{John2017}, we prescribe a (ramped up) pulsating horizontal inflow profile with period $T_p=0.917$ s:
\begin{align*}
u(y,t)|_{x=0} = \frac12\Bigg[1-\Big(\frac{2y}{H}-1\Big)^2\Bigg]U(t) f(t)\, ,
\end{align*}
with $U(t)$ as in \eqref{ramp}, $t^{\star}=T_p$. The function $f(t)$ is a trigonometric polynomial fitted from internal carotid artery (ICA) measurements \cite{Thiriet2008}:
\begin{align*}
f(t) = 10^{-4}\Bigg[\frac{a_0}{2} - \sum_{k=1}^{12}a_k\cos\left(\frac{2k\pi t}{T_p}\right) + b_k\sin\left(\frac{2k\pi t}{T_p}\right)\Bigg]\, ,
\end{align*}
where $a_0 = 13349$, and the remaining coefficients are listed on Table \ref{FourierData}. That leads to a peak inflow speed of $U_{\mathrm{peak}} = 0.5$ m/s. At the outlet ($x=6H$), we prescribe the outflow condition
\begin{align}
    (\nu\nabla\ve{u})\ve{n} - p\ve{n} = \ve{0}\, ,
\label{pseudotraction}
\end{align}
while the bottom and top walls are no-slip boundaries ($\ve{u}=\ve{0}$).
\begin{table}[ht!]
\caption{Trigonometric coefficients for the ICA inflow waveform ($a_0 = 13349$).}
{\scriptsize
\begin{tabular}{c|cccccccccccccc}
   $k$  & 1 &  2 & 3 & 4 & 5 & 6 & 7 & 8 & 9 & 10 & 11 & 12 \\
   \hline
   $a_k$  & 695.98 & 905.24 & 452.42 & 754.46 & 164.6 & 133.57 & 79.4 & 74.40 & 36.72 & 38.80 & 42.34 & 3.67  \\
  \hline
     $b_k$ & $-2067.65$  & $-496.92$  & $-340.74$ & 202.96 & 98.82 & 249.55 & 81.7 & 67.84 & 9.95 & 15.34 & $-3.49$ & $-3.22$\\
\hline
 \end{tabular}}
 \label{FourierData}
\end{table}

\begin{figure}[ht!]
\centering
\includegraphics[trim = 0 0 0 0,clip, width = .8\textwidth]{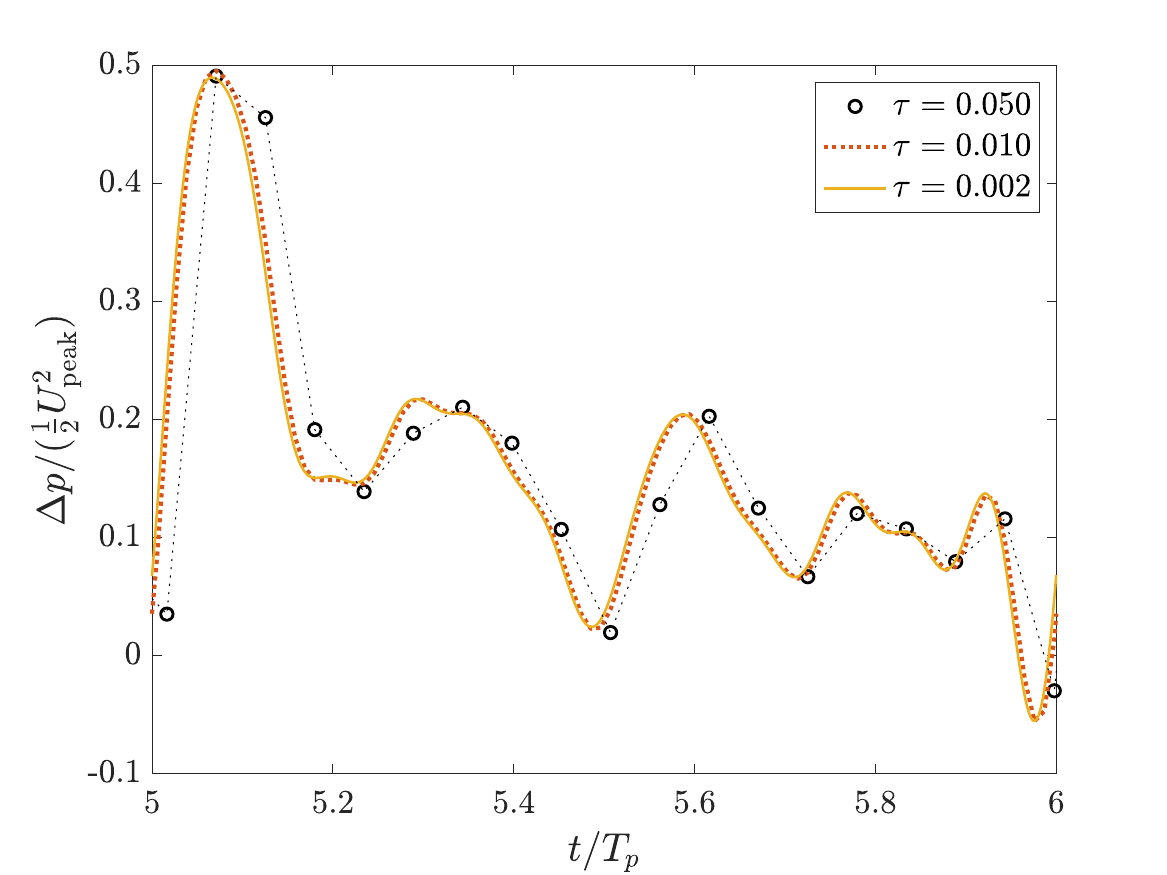}
\caption{Carreau fluid through Idealised ICA aneurysm: normalised pressure difference (inlet minus outlet) over one cardiac cycle, for different time-step sizes $\tau$. All three computations are based on the monolithic IMEX scheme in GL form, which allows enforcing the outflow condition \eqref{pseudotraction} naturally.}
\label{aneurysm_dp}
\end{figure}

\begin{figure}[ht!]
\centering
\includegraphics[trim = 160 270 130 250,clip, width = .85\textwidth]{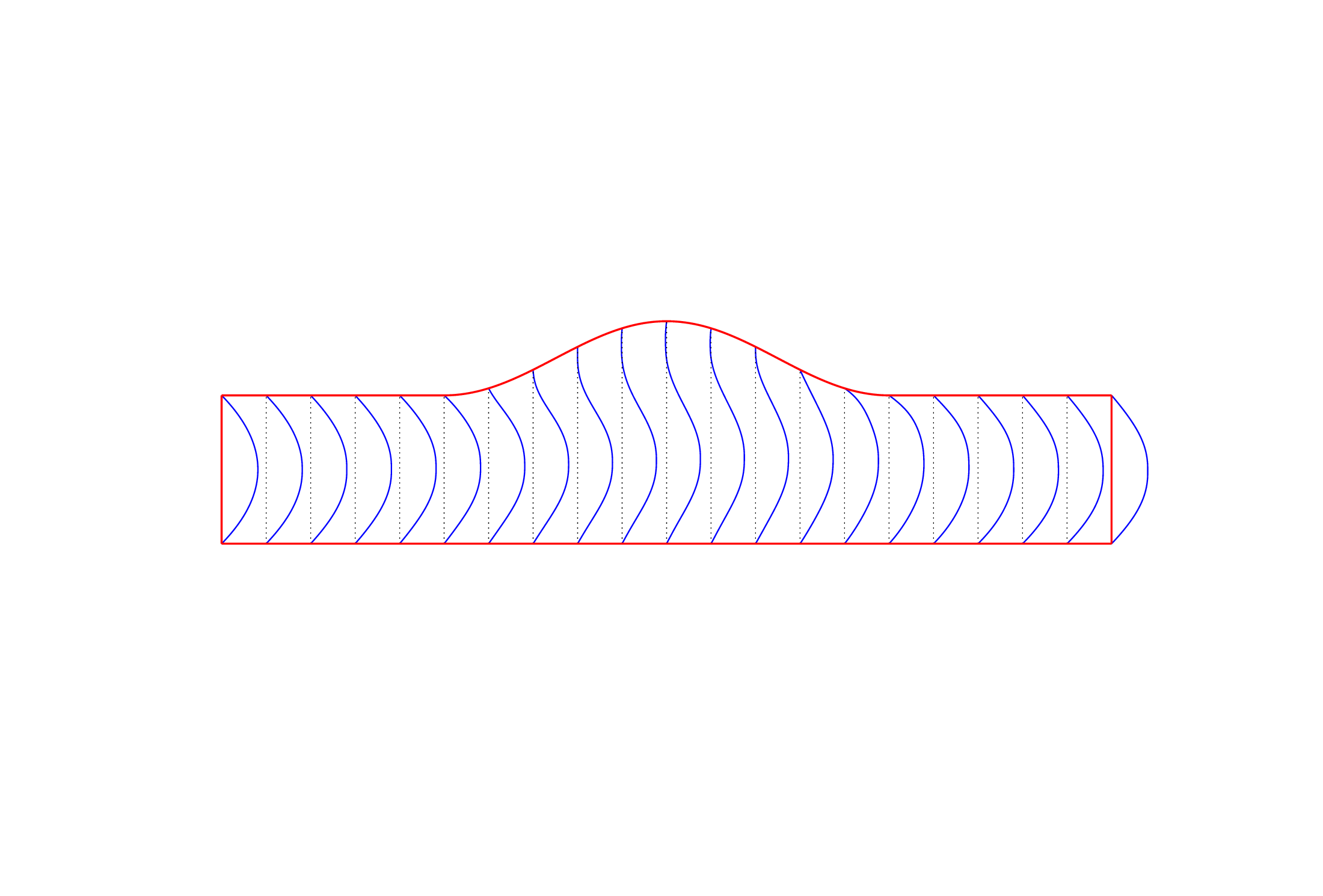}
\caption{Carreau fluid through Idealised ICA aneurysm: horizontal velocity profiles at $t=6T$, for a time-step size $\tau=0.002$.}
\label{aneurysmProfiles}
\end{figure}

For this test we consider the GL formulation, as it offers a natural way to enforce the outflow condition \eqref{pseudotraction}. We use the monolithic IMEX-GL scheme \eqref{decoupledGL} with different time-step sizes $\tau$ and a mesh with 38,400 elements. Figure \ref{aneurysm_dp} shows, over one cardiac cycle, the mean pressure drop across the aneurysm:
\begin{align*}
    \Delta p = \frac{1}{H}\left(\int_{\Gamma_{\mathrm{in}}} p \, \mathrm{d}\Gamma - \int_{\Gamma_{\mathrm{out}}} p \, \mathrm{d}\Gamma\right),
\end{align*}
which is an important clinical marker \cite{Pacheco2022IJNMBE}. The results confirm very good temporal stability for the numerical solutions, even for time-step sizes as large as $\tau = 0.05$. For $t=6T$, Figure \ref{aneurysmProfiles} shows the horizontal velocity profiles at different sections of the vessel, showcasing also good spatial stability. 

\section{Concluding remarks}\label{sec_Conclusion}
In this work, we have presented, analysed and tested different implicit-explicit schemes for the variable-viscosity Navier--Stokes equations. While IMEX methods for constant viscosity focus on making  the convective velocity (or different stabilisation terms associated to convection) explicit, we split also the viscous term into implicit and explicit parts. This renders the velocity subsystem block-diagonal, thereby enabling simpler, more efficient linear solvers. The temporal stability of the resulting schemes varies dramatically according to the term made explicit, and in which form. To investigate that, we have carried out a semi-discrete analysis of prototypical first-order  schemes. Our results show that, when considering a stress-divergence (SD) description of the viscous term, naive extrapolations can cause instabilities. While unconditionally stable IMEX-SD schemes are possible, we have shown that they require a pre-computed viscosity field, which is not realistic in many applications. On the other hand, a generalised Laplacian (GL) formulation of the viscous term leads to stable IMEX methods requiring only known viscosity values. For a monolithic and a fractional-step IMEX-GL scheme, respectively, we have proved unconditional and conditional stability. Our numerical tests, however, suggest that the time-step restriction derived for the fractional-step variant may be an analytical artefact.

Our theoretical and numerical results contribute towards breaking the paradigm that considers \textit{inevitable} the use of fully coupled, elasticity-like formulations for variable-viscosity flows. The GL variants presented in this work have provided  accurate and stable results, performing well even for problems with non-smooth solutions. Several questions remain open though. Proposing and analysing second- and higher-order schemes is a challenging and interesting open question. In addition, analysing the stability and convergence of fully discrete versions of the schemes proposed in this work is also of interest, especially for the case of nonlinear viscosity models.


\section*{Acknowledgments}
 The work of GRB has been partially supported by the Leverhulme Trust Research Project Grant No. RPG-2021-238.  EC acknowledges the support given by the Agencia Nacional de Investigaci\'on y Desarrollo (ANID) through the project FONDECYT 1210156.

\bibliographystyle{unsrtnat}
\bibliography{references}%

\end{document}